\documentclass[11pt,reqno]{amsart}
\usepackage{amsmath, amsthm, amssymb, amscd, amsfonts}
\usepackage[all]{xy}
\usepackage{color}

\def\pf{\begin{proof}}
\def\epf{\end{proof}}

\newcommand{\ev}{\rm ev}

\newcommand\id{\operatorname{id}}
\newcommand{\Ima}{\mbox{\rm Im\,}}

\newcommand{\ku}{{\Bbbk}}
\newcommand{\Z}{{\mathbb Z}}
\newcommand{\N}{{\mathbb N}}
\newcommand{\R}{{\mathbb R}}
\newcommand{\C}{{\mathbb C}}
\newcommand{\Q}{{\mathbb Q}}
\newcommand{\X}{\left(\frac c{cz+d}\right)}
\newcommand{\ie}{i.e.\ }
\newcommand{\trsf}[3]{(#1{\vert}_{#2}#3)}%..........(f|M)

\theoremstyle{plain}
\newtheorem{theo}{Theorem}[subsection]
\newtheorem{lemma}[theo]{Lemma}
\newtheorem{coro}[theo]{Corollary}
\newtheorem{prop}[theo]{Proposition}

\theoremstyle{definition}
\newtheorem{defi}[theo]{Definition}
\newtheorem{defis}[theo]{Definitions}
\newtheorem{rem}[theo]{Remark}
\newtheorem{rems}[theo]{Remarks}

\newtheorem{exas}[theo]{Examples}

\newtheorem{problem}[theo]{Problem}

\newtheorem{nota}[theo]{Notations}

\begin{document}
\title[Invariants of formal pseudodifferential operator algebras]{\bf Invariants of formal pseudodifferential operator algebras and algebraic modular forms}
\author{Fran\c cois Dumas}
\address[F.~Dumas]{Universit\'e Clermont Auvergne, Laboratoire de Math\'ematiques Blaise Pascal (UMR 6620 du CNRS), F-63178, Aubi\`ere (France)}
\email{Francois.Dumas@uca.fr} 
\author{Fran\c cois Martin}
\address[F.~Martin]{Universit\'e Clermont Auvergne, Laboratoire de Math\'ematiques Blaise Pascal (UMR 6620 du CNRS), F-63178, Aubi\`ere (France)}
\email{Francois.Martin@uca.fr}

\begin{abstract}
We study from an algebraic point of view the question of extending an action of a group \(\Gamma\) on a commutative domain \(R\) to a formal pseudodifferential operator ring 
\(B=R(\!(x\,;\,d)\!)\) with coefficients in \(R\), as well as to some canonical quadratic extension \(C=R(\!(x^{1/2}\,;\,\frac 12 d)\!)_2\) of \(B\). 
We give a necessary and sufficient condition of compatibility between the action and the derivation $d$  of $R$ for such an extension to exist, and we determine all possible extensions of the action to \(B\) and \(C\).
We describe under suitable assumptions the invariant subalgebras \(B^\Gamma\) and \(C^\Gamma\) as Laurent series rings with coefficients in \(R^\Gamma\). 

The main results of this general study are applied in a numbertheoretical context to the case where
\(\Gamma\) is a subgroup of \({\rm SL}(2,\C)\) acting by homographies on an algebra \(R\) of functions in one complex variable.
Denoting by \(M_j\) the vector space of algebraic modular forms in $R$ of weight \(j\) (even or odd),
we build for any nonnegative integer \(k\) a linear isomorphism between the subspace \(C_k^\Gamma\) of invariant operators of order \(\geq k\) in \(C^\Gamma\)
and the product space \(\mathcal{M}_k=\prod_{j\geq k}M_j\), which can be identified with a space of algebraic Jacobi forms of weight \(k\).
It results in particular a structure of noncommutative algebra on \(\mathcal M_0\) 
and an  algebra isomorphism \(\Psi:\mathcal M_0\to C_0^\Gamma\), whose restriction to the particular case of even weights
was previously known in the litterature. We study properties of this correspondence combining arithmetical arguments and the use of the algebraic results of the first part of the article.\end{abstract}

\date\today
\thanks{The authors are partially funded by the project CAP 20-25, I-SITE Clermont.}
\maketitle

\section*{Introduction}
Several studies in deformation or quantization theories have shown in various contexts that significant combinations
of Rankin-Cohen brackets on modular forms of even weights correspond by isomorphic transfer to the noncommutative composition product in
some associative algebras of invariant operators (see for instance \cite{Z2}, \cite{CMZ}, \cite{UU},  \cite{BXY}, \cite{CL}, \cite{MP}, \cite{Yao2}, \cite{lee2}). 
The main goal of this paper is to produce such a correspondence 
in a formal algebraic setting for modular forms of any weights (even or odd), which allows the construction to be applied to Jacobi forms.
\medskip

The first part deals with the general algebraic problem to extend a group action from a ring $R$ 
to a ring of formal pseudodifferential operators with coefficients in $R$, and to describe the subring of invariant operators. 
More precisely, let $R$ be a commutative domain of characteristic zero containing the subfield of rational numbers, $U(R)$ 
its group of invertible elements, $d$ a derivation of $R$ and
$\Gamma$ a group acting by automorphisms on $R$. We denote by $B=R(\!(x\,;\,d)\!)$ the ring of formal 
pseudodifferential operators with coefficients in $R$. The elements of $B$ are the Laurent power series in one
indeterminate $x$ with coefficients in $R$ and the noncommutative product in $B$
is defined from the commutation law: 
\begin{center}
$xf=\sum_{i\geq 0}d^i(f)x^{i+1}$ for any $f\in R$.\end{center}
We also consider some quadratic extension 
$C=R(\!(y\,;\,\delta)\!)_2$ of $B$, with \(y^2=x\), \(\delta=\frac12d\) and commutation law:\begin{center}
$\textstyle{yf=\sum_{k\geq 0}\left({\textstyle{{\frac{(2k)!}{2^k(k!)^2}}}}\right)\delta^k(f)y^{2k+1}}$ for any $f\in R$.\end{center}
This type of skew power series rings was already introduced  in various ringtheoretical papers, see references in \cite{FD}. 
We prove (theorem \ref{extensionB}) that the action of $\Gamma$
on $R$ extends to an action by automorphisms on $B$ if and only if there exists a multiplicative
1-cocycle $p:\Gamma\to U(R)$ such that $\gamma d\gamma^{-1}=p_\gamma.d$ for any $\gamma\in\Gamma$.
We describe under this assumption all possible extensions of the action; they are parametrized by the arbitrary choice of a map
$r:\Gamma\to R$ satisfying some compatibility condition related to $p$. We give (theorem \ref{casimportantB}) sufficient conditions for the
ring $B^\Gamma$ of invariant operators to be described as a ring of formal pseudodifferential operators with coefficients in $R^\Gamma$.
We prove (theorems \ref{extensionC} and \ref{casimportantC}) similar results for the ring $C$; in this case a necessary and sufficient condition to extend the action of $\Gamma$ from $R$ to $C$ 
is the existence of some multiplicative 1-cocycle $s:\Gamma\to U(R)$ satisfying $\gamma d\gamma^{-1}=s^2_\gamma.d$ for any $\gamma\in\Gamma$. 
The subring $B$ is then necessarily stable under the extended actions, and \(B^\Gamma\) is a subring of \(C^{\Gamma}\).\medskip

These general results are applied in the second part of the paper 
to the case of the complex homographic action. We fix a $\C$-algebra $R$ of functions in one variable $z$ 
stable under the standard derivation $\partial_z$ and a subgroup $\Gamma$ of ${\rm SL}(2,\C)$ acting on $R$ by:\begin{center} 
$(f.\gamma)(z)=f(\frac{az+b}{cz+d})$ for any $f\in R$ and $\gamma=\left(\begin{smallmatrix}a &b\\ c&d \end{smallmatrix}\right)\in\Gamma$.\end{center}
Applying the previous general construction for the derivation $d=-\partial_z$ and the 1-cocycles 
$s~:~\gamma\mapsto cz+d$ and $p=s^2$, we define an action of $\Gamma$ on 
$C=R(\!(y\,;\,-\frac12\partial_z)\!)_2$. 
We give in theorems \ref{extensionBC} and \ref{extensionyk}  a combinatorial description of this action, whose restriction to $B=R(\!(x\,;\,-\partial_z)\!)$ with \(x=y^2\) corresponds to the action
already introduced in section 1 of \cite{CMZ}. 
We consider for any $k\in\Z$ the subspace $C_k^{\Gamma}$ of elements in $C^{\Gamma}$ whose valuation related to $y$ is greater or equal to $k$.
It is easy to check that the image of the canonical projection $\pi_k:C_k^{\Gamma}\to R$ is included in the subspace $M_{k}$ of algebraic modular
forms of weight $k$. We construct a family \(\!(\psi_k)_{k\geq 0}\) of splitting maps  $\psi_k:M_{k}\to C_k^\Gamma$
which gives rise canonically (theorem \ref{*positif}) to a vector space isomorphism:\medskip

\centerline{$\Psi: \mathcal M_{0}\to C_0^{\Gamma}$,\ \ with $\mathcal M_{0}=\prod_{k\geq 0}M_{k}$.}\medskip

This correspondence makes it possible on the one hand to identify $C_k^{\Gamma}$ with a vector space of algebraic Jacobi forms of weight $k$ (see corollary \ref{jacobi2}), 
on the other hand to obtain by transfer a noncommutative product on  $\mathcal M_{0}$. This multiplication on modular forms of nonnegative weights (even or odd) can be described in terms of linear combinations of Rankin-Cohen brackets with combinatorial rational coefficients:\begin{center}
$f\star g=\Psi^{-1}(\psi_k(f).\psi_\ell(g)\!)=\sum_{n\geq 0}\alpha_n(k,\ell)[f,g]_n$\end{center}
for any $f\in M_k$, $g\in M_\ell$, $k,\ell\geq 0$.
The restriction to modular forms of even nonnegative weights gives the isomorphism:\medskip

\centerline{$\Psi_2: \mathcal M^{\ev}_0\to B_0^{\Gamma}$,\ \ with $\mathcal M^{\ev}_0=\prod_{k\geq 0}M_{2k}$ and \(B_0^\Gamma=B\cap C_0^\Gamma\),}\medskip

considered in section 3 of \cite{CMZ} (see also \cite{UU}, \cite{Yao}).
We discuss  at the end of paragraph \ref{23} the obstructions for possible extensions of theorem \ref{*positif} to modular forms of negative weights.
The purpose of paragraph \ref{24} is to go further in the study of the isomorphism $\Psi_2$ using the structure of $B_0^{\Gamma}$
deduced from the general theorem \ref{casimportantB}. More precisely, assuming that $R$ contains an invertible modular form \(\chi\) of weight 2,
\(B_0^\Gamma\) is the skew power series algebra \(R^{\Gamma}[[u\,;\,D]]\) with coefficients in \(R^\Gamma\) where \(u\) denotes the invariant operator \(x\chi\)
and \(D\) is the derivation \(-\chi^{-1}\partial_z\) of \(R^{\Gamma}\). A natural question is then to describe as a sequence of modular forms the inverse image by \(\Psi_2\)
of any invariant operator in \(B_0^\Gamma\), or equivalently of any power \(u^k\) for \(k\geq 0\).
The value of $\Psi_2^{-1}(u^k)$ is entirely determined by some family \(\!(g_{k,k+i})_{i,k\geq 0}\) of modular forms of weight \(2k+2i\), which are zero if \(i\) is odd,
and which can be calculated for even \(i\)'s in terms of Rankin-Cohen brackets of powers of the modular form \(\chi\). Paragraph \ref{25} explores some ways to extend the previous results to negative odd weights.
We finally mention in paragraph \ref{26} as an open question the possibility of considering the whole study of the second part of this paper for other choices of the parameter $r$ in the  extension
of the initial action from functions to operators.

\section{Actions and invariants on pseudodifferential
and quadratic pseudodifferential operator
rings}\label{firstsection}

\subsection{Definitions and preliminary results}

We fix a commutative domain $R$ of characteristic zero containing $\Q$. We denote by $U(R)$ the group of invertible elements in $R$.
We denote by $R\langle\!\langle x\rangle\!\rangle$ the left $R$-module 
of formal power series $\sum_{i\geq 0}f_ix^i$ with $f_i\in R$ for
any $i\geq 0$. It is well known that, for any nonzero derivation $d$ of $R$,
we can define an associative noncommutative product on $R\langle\!\langle x\rangle\!\rangle$
from the commutation law:
\begin{equation}\label{PDOx1}
xf=fx+d(f)x^2+d^2(f)x^3+\cdots=\sum_{n \geq 0} d^n(f) x^{n+1}
\hspace{0.3cm} \text{ for any } f \in R,\end{equation} and more generally for any integer $i\geq 0$:
\begin{equation}\label{PDOxi}
x^{i+1}f=fx^{i+1}+(i+1)d(f)x^{i+2}+\textstyle{\frac{(i+1)(i+2)}2}d^2(f)x^{i+3}+\cdots=\sum_{n
\geq i} {\textstyle{\binom ni}} d^{n-i}(f) x^{n+1}.
\end{equation}
We obtain so a noncommutative ring denoted by $B_0=R[[x; d]]$.
The element $x$ generates a two-sided ideal in $B_0$ and
we can consider the localized ring $B=R(\!(x;d)\!)$ of $B_0$ with
respect of the powers of $x$. The elements of $B$ are the Laurent
series $q=\sum_{i\geq m}f_ix^i$ with $m\in\Z$ and $f_i\in R$ for
any $i\geq m$. The valuation $v_x(q)$ of $q$ is  defined by
$v_x(q)=\min\{i\in \Z\,;\, f_i\not=0\}$ for $q\not=0$ and $v_x(0)=\infty$. The
map $v_x$ is a valuation function $B\to\Z \cup \{\infty\}$ and $B_0$ is the subring $\{q\in B\,;\,
v_x(q)\geq 0\}$. It follows in particular that $B_0$ and $B$ are
domains. The invertible elements of $B$ are the series $q=\sum_{i \geq m} f_i x^i$ such that $f_m \in U(R)$. \medskip

\begin{defi}\label{PDO} The noncommutative domain $B=R(\!(x;d)\!)$ is called the ring
of formal pseudodifferential operators in $d$ with coefficients in $R$.
\end{defi}

Denoting by $\partial$ the derivation $-d$, relation
\eqref{PDOx1} can be rewritten:
\begin{equation}\label{schur3}
x^{-1}f=fx^{-1}+\partial (f) \ \text{ for any } f\in R.
\end{equation}
It is well known that this relation defines an associative noncommutative product in the subset $A$ of
polynomials in the variable $t=x^{-1}$ with coefficients in $R$
(see for instance \cite{GW} pages 11 and 19). This subset $A$ is then a subring of $B$, denoted by $A=R[t\,;\,\partial]$.
\begin{defi}\label{DO} The noncommutative domain $R[t\,;\,\partial]$ is called the ring of formal differential operators in $\partial$ with coefficients in $R$. 
\end{defi}

Denoting by $y$ instead of $x$ the variable in the left $R$-module $R\langle\!\langle y\rangle\!\rangle$,
another structure of noncommutative ring (see example 1.3.(d) in \cite{FD}) can be defined on $R\langle\!\langle y\rangle\!\rangle$ for any nonzero
derivation $\delta$ of $R$ from the commutation law:
\begin{equation}\label{QPDOy1}
{\textstyle{yf=fy+\delta(f)y^3+\frac32\delta^2(f)y^5+
\frac52\delta^3(f)y^7+\frac{35}8\delta^4(f)y^9+\cdots}}=\sum_{k
\geq 0}\left({\textstyle{{\frac{(2k)!}{2^k(k!)^2}}}}\right)\delta^k(f)
y^{2k+1},\end{equation} and more generally  for any integer $i\geq 0$:
\begin{equation}\label{QPDOyi}
y^{i+1}f=fy^{i+1}+\sum_{k \geq 1}\left(\prod_{j=1}^k\frac{2(j-1)+i+1}{j}\right)\delta^k(f)y^{2k+i+1}  \ \text{ for any } f\in R.
\end{equation}
This ring is denoted by $C_0=R[[y;\delta]]_2$. As in the previous case, $C_0$ can be
embedded in the localization $C=R(\!(y;\delta)\!)_2$ of Laurent series $q=\sum_{i\geq m}f_iy^i$ with $m\in\Z$ and $f_i\in R$ for
any $i\geq m$. The valuation $v_y(q)$ is  defined by
$v_y(q)=\min\{i\in \Z\,;\, f_i\not=0\}$ for $q\not=0$ and $v_y(0)=\infty$. The
map $v_y$ is a valuation function $C\to\Z \cup \{\infty\}$ and $C_0$ is the subring $\{q\in B\,;\,
v_y(q)\geq 0\}$. In particular $C_0$ and $C$ are domains. The invertible elements of $C$ are the series $q=\sum_{i \geq m} f_i y^i$ such that $f_m \in U(R)$.\medskip

\begin{defi}\label{QPDO} The noncommutative domain $C=R(\!(y;\delta)\!)_2$ is called the ring
of quadratic formal pseudodifferential operators in $\delta$ with coefficients in $R$.
\end{defi}

By elementary calculations it follows from \eqref{QPDOyi} that:
\begin{equation}\label{QPDOy2}
y^2f=fy^2+2\delta(f)y^4+4\delta^2(f)y^6+\cdots=\sum_{k \geq
1}(2\delta)^{k-1}(f)y^{2k} \ \text{ for any } f\in R,
\end{equation}
which is equivalent to
\begin{equation}\label{QPDOy-2}
y^{-2}f=fy^{-2}-2\delta(f) \ \text{ for any } f\in R.
\end{equation}
This observation gives rise to the following result.
\begin{prop}\label{BC}Let $R$ be a commutative domain of characteristic zero containing $\Q$.
Let $d$ be a nonzero derivation of $R$. 
\begin{itemize}
\item[{\rm (i)}] The ring of quadratic pseudodifferential operators $C=R(\!(y\,;\,\delta)\!)_2$ 
for  $\delta=\frac12d$ contains as a subring the ring of pseudodifferential operators $B=R(\!(x\,;\,d)\!)$
for $x=y^2$, and therefore the ring of
differential polynomials $A=R[t\,;\,\partial]$ for $t=x^{-1}=y^{-2}$ and
$\partial=-d$.\begin{equation*}
A=R[t\,;\,-d]\ \subset \ B=R(\!(x\,;\,d)\!)\ \subset \
C=R(\!(y\,;\, \textstyle\frac12d)\!)_2.\end{equation*}
\item[{\rm (ii)}] We have $v_y(q)=2v_x(q)$ for any $q\in B$, and the subrings
$B_0=R[[x\,;\,d]]$
and $C_0=R[[y\,;\,\delta]]_2$ 
satisfy $B_0=B\cap C_0$. 
\medskip
\item[{\rm (iii)}] 
We have $C=B\oplus By$,  with  $By=yB$, 
$(By)(By)=B$ and $(By)B=B(By)=By$.
\end{itemize}\end{prop}
{\it Proof.} Follows by obvious calculations from identities \eqref{QPDOy2} and \eqref{QPDOy-2}.\qed
\bigskip

The following proposition shows that
any automorphism of $B$ (respectively $C$) stabilizing $R$ is continuous with respect of the $x$-adic
(respectively the $y$-adic) topology. A similar property is proved in \cite{AD1} for different
commutation laws and in the particular case where $R$ a field. We
start with a preliminary technical result.

\begin{lemma}\label{root} Data and notations are those of the previous proposition.
Let $\ell$ be a positive integer. \begin{itemize} \item[{\rm(i)}]
 For any $q\in B$ of the form $q=1+\sum_{i\geq 1}f_ix^i$
 with $f_i\in R$ for any $i\geq 1$, there exists an element
 $r$ in $B$ such that $q=r^{\ell}$.  \item[{\rm(ii)}] For any $q\in C$ of the form
 $q=1+\sum_{i\geq 1}f_iy^i$
 with $f_i\in R$ for any $i\geq 1$, there exists an element
$r$ in $ C$ such that $q=r^{\ell}$.
\end{itemize}\end{lemma}

{\it Proof.} The proof is the same in the two cases and we write it for $B=R(\!(x;d)\!)$. 
We consider in $B$ an element of
nonnegative valuation $r=\sum_{i\geq 0} g_ix^i$ with $g_i\in R$
for any $i\geq 0$. For any integer $\ell\geq 1$, we denote
$r^{\ell}=\sum_{i\geq 0} g_{\ell,i}x^i$ with $g_{\ell,i}\in R$ for
any $i\geq 0$. By a straightforward induction using
\eqref{PDOxi}, we check that $g_{\ell,i}=\ell
g_0^{\ell-1}g_i+h_{\ell,i}$ where the rest $h_{\ell,i}$
depends only on previous $g_0,g_1,\ldots,g_{i-1}$ and their images by
some powers of $d$. Then for any sequence $(f_i)_{i\geq 1}$ of elements in
$R$, we can determine inductively a unique sequence $(g_i)_{i\geq 0}$
such that $g_0=1$ and $g_{\ell,i}=f_i$ for any $i\geq 1$, that is
$(\sum_{i\geq 0} g_ix^i)^{\ell}=1+\sum_{i\geq 1}f_ix^i$. \qed

\begin{prop}\label{continuity}Data and notations are those of the previous proposition.
\begin{itemize} \item[{\rm(i)}]
Let $\gamma$ be an automorphism of $B$ whose restriction to $R$ is
an automorphism of $R$. Then $v_x(\gamma(q))=v_x(q)$ for any $q\in
B$. In particular the restriction of $\gamma$ to $B_0$
determines an automorphism of $B_0$.\item[{\rm(ii)}] Let $\gamma$ be an automorphism of $C$
whose restriction to $R$ is an automorphism of $R$. Then
$v_y(\gamma(q))=v_y(q)$ for any $q\in C$. \item[{\rm(iii)}] Let
$\gamma$ be an automorphism of $C$ whose restriction to $R$ is an
automorphism of $R$. Then the restriction of $\gamma$ to $B$
determines an automorphism of $B$.
\end{itemize}
\end{prop}

{\it Proof.} Let $\gamma$ be an automorphism of $B$ such that
$\gamma(R)= R$. This implies $\gamma^{-1}(R)=R$. We introduce
$m=v_x(\gamma(x)) \in \Z$ and suppose that $m<0$. Then the element
$z=1+x^{-1}$ satisfies $\gamma(z) = 1+\gamma(x)^{-1} \in B_0$. We
set an integer $\ell\geq 2$ and consider by point (i) of lemma
\ref{root} an element $r \in B_0$ such that
$r^\ell=\gamma(z)$. Applying the automorphism $\gamma^{-1}$ we have
$v_x(z) = \ell v_x(\gamma^{-1}(r))$, which gives a contradiction because
$v_x(z)=-1$ by definition. Hence we have proved that $m \geq 0$.
It follows in particular that $\gamma(B_0)\subset B_0$. By the
same argument for $\gamma^{-1}$, we conclude that $\gamma(B_0)=
B_0$. In other words, the restrictions of $\gamma$ and $\gamma^{-1}$ to $B_0$ determine
automorphisms of $B_0$. We prove now that $m=1$. In $B_0$, we can
write $\gamma(x) = f (1+w) x^m$, with $f \in U(R)$ (because $\gamma(x)$ is invertible in $B$) and $w \in B_0$
such that $v_x(w) \geq 1$. It follows that $x=\gamma^{-1}(f)
\gamma^{-1}(1+w) \gamma^{-1}(x)^m$, and then $v_x(\gamma^{-1}(f))
+ v_x(\gamma^{-1}(1+w)) + mv_x(\gamma^{-1}(x)) = 1$. Since
$\gamma^{-1}(R)=R$, we have $v_x(\gamma^{-1}(f))=0$. The element $1+w$ is invertible in $B_0$, hence
$\gamma^{-1}(1+w)$ is invertible in $B_0$, therefore $v_x(\gamma^{-1}(1+w)) = 0$. We conclude that
$mv_x(\gamma^{-1}(x))=1$, then $m=1$ and the proof of point (i)
is complete. The proof of point (ii) is similar replacing $B_0$ by
$C_0$, $x$ by $y$ and $v_x$ by $v_y$.\medskip

Let $\gamma$ be an automorphism of $C$ such that $\gamma(R)=R$. We
deduce from relation \eqref{QPDOy-2} that
$\gamma(y)^{-2}\gamma(f)-\gamma(f)\gamma(y)^{-2}=-\gamma(d(f))$
for any $f\in R$. With notation $z=\gamma(y)^{-2}$, we obtain 
$zf-fz=-\gamma d\gamma^{-1}(f) \in R$ for any $f\in R$. By point
(ii), we have $v_y(z)=-2$. Using point (iii) of proposition \ref{BC}, we
have $z=q+q'y$ with $q,q'\in B$, $v_y(q)=-2$, $v_y(q')\geq
-2$. Since $qf-fq\in B$ for any $f\in R$, we deduce from previous
identity that $q'yf-fq'y\in B$ for any $f\in R$. Suppose that
$q'\not=0$ and denote $q'y=\sum_{i\geq \ell}f_{2i+1}y^{2i+1}$, with
$\ell\geq -1$, $f_{2i+1}\in R$ for any $i\geq -1$ and
$f_{2\ell+1}\not=0$. Using \eqref{QPDOyi}, we have in $C$ for
any $f\in R$ the development: $q'yf-fq'y=(2\ell+1)f_{2\ell+1}\delta(f)y^{2\ell+3}+\cdots$ which is incompatible with the fact that $q'yf-fq'y \in B$ and $\delta$ is a nonzero derivation.
Then we
have necessarily  $q'=0$, that is $z\in B$. In other words,
$\gamma(x^{-1})\in B$. Hence $\gamma(x)\in B$ and the proof is
complete. \qed

%%%%%%%%%%%%%%%%%%%%%%%%%%%%%%%%%%%%%%%%%%%%%%%%%%%%%%%%%%%%%%%%%%%%%%%%%%%%%%%%%%%%%%%%%%%%%%%%%%%%%%%%%%%%%%%%%%%%%%%%%%%%%%%%%%%ù

\subsection{Extension to $B$ and $C$ of actions by automorphisms on $R$}\label{12}

\begin{nota}\label{mainnota}
We take all data and notations of proposition \ref{BC}.
We consider a group $\Gamma$ acting by automorphisms on the ring $R$. 
We denote this action on the right:
\begin{equation}\label{actionautom}
(f\cdot\gamma)\cdot\gamma'=f\cdot \gamma\gamma' \quad\text{for all} \ f\in R,\ \gamma,\gamma'\in\Gamma.\end{equation}
A 1-cocycle for the action of $\Gamma$ on the group
$U(R)$ is a map $s: \Gamma\to U(R)\,,\, \gamma\mapsto s_{\gamma}$ satisfying:
\begin{equation}\label{cocycle}s_{\gamma\gamma'}=(s_\gamma\cdot\gamma')s_{\gamma'}
\quad \text{ for all } \gamma,\gamma'\in\Gamma.\end{equation} We denote by $Z^1(\Gamma,U(R))$
the multiplicative abelian group of such 1-cocycles. \medskip

We answer the following two questions: give necessary and sufficient conditions for the existence
of an action of $\Gamma$ by automorphisms on $B$ or $C$ extending the given action on $R$, 
and describe all possible extended actions. We need some definitions.
\end{nota}

\begin{defis}\label{troisdef}Data and notations are those of \ref{mainnota}.
\begin{itemize}
\item[(i)] The action of $\Gamma$ on $R$
is said to be $d$-compatible when for any $\gamma\in\Gamma$, there exists
$p_{\gamma}\in U(R)$ such that:
\begin{equation}\label{dcomp}
 d(f)\cdot\gamma=p_\gamma d(f\cdot\gamma) \ \ \text{ for any } f\in R. 
\end{equation} This condition 
defines uniquely a map $p:\Gamma\to U(R)\,,\,\gamma\mapsto p_{\gamma}$, and 
$p\in Z^1(\Gamma,U(R))$. This map $p$ is called the $1$-cocycle associated to the $d$-compatibility.
\item[(ii)] The action of $\Gamma$ on $R$ is said to be quadratically $d$-compatible when 
there exists an element $s\in Z^1(\Gamma,U(R))$ such that:
\begin{equation}\label{quaddcomp}d(f)\cdot\gamma=s_\gamma^2 d(f\cdot\gamma) \ \ 
\text{for any } \gamma\in\Gamma \text{ and any } f\in R.\end{equation}
Such a map $s$ is called a $1$-cocycle associated to the quadratic $d$-compatibility.
It is not necessarily unique since any
map $s'=\epsilon s$ where $\epsilon$ is a multiplicative
function $\Gamma~\to~\{-1,+1\}$ is another element of
$Z^1(\Gamma,U(R))$ satisfying \eqref{quaddcomp}.\end{itemize}
\end{defis}

\begin{rem}
Any quadratically $d$-compatible action is also $d$-compatible.
Since $\delta=\frac12d$ and $\partial=-d$, the $d$-compatibility is obviously equivalent to 
the $\delta$-compatibility or the $\partial$-compatibility. 
\end{rem}

\begin{exas}\label{exemplesactions}
\begin{itemize}
\item[1.] For any commutative field $\ku$
of characteristic zero, the group $\Gamma=\ku\rtimes~\ku^{\times}$
for the product $(\mu,\lambda)(\mu',\lambda')=(\lambda\mu'+\mu,\lambda\lambda')$
acts by $\ku$-automorphisms on the domain $R=\ku[z]$ by $(f\cdot\gamma)(z)=f(\lambda z + \mu)$ for $\gamma=(\mu,\lambda)\in \Gamma$. 
It is clear that
$\partial_z(f)\cdot\gamma=\lambda^{-1}\partial_z(f\cdot\gamma)$ for any polynomial $f\in R$. Hence this action is
$\partial_z$-compatible, with $p\in Z^1(\Gamma,\ku^\times)$ defined by
$p:(\mu,\lambda)\mapsto\lambda^{-1}$. If $\ku$ is not algebraically
closed, the action is not necessarily quadratically
$\partial_z$-compatible. \item[2.]Let
$\ku$ be any commutative field of characteristic zero and
$R=\ku(z)$ the field of rational functions in one variable with
coefficients in $\ku$. The group $\Gamma={\rm SL}(2,\ku)$ acts by $\ku$-automorphisms 
on $R$ by $(f\cdot\gamma)(z)=f(\frac{\lambda z + \mu}{\eta z + \xi})$
where $\gamma=(\begin{smallmatrix} \lambda &\mu\\ \eta
&\xi\end{smallmatrix})\in\Gamma$. We have $\partial_z(f)\cdot\gamma=(\eta z+\xi)^2\partial_z(f\cdot\gamma)$ for any rational function $f\in R$.
Hence this action is quadratically $\partial_z$-compatible, with $s\in
Z^1(\Gamma,R^\times)$ defined by $s:(\begin{smallmatrix} \lambda
&\mu\\ \eta &\xi\end{smallmatrix})\mapsto\eta z+\xi$. This
kind of action is the main object of the second part of the paper.\end{itemize}
\end{exas}

We are now able to answer for $B$ to the questions formulated at
the beginning of this paragraph.

\begin{theo}\label{extensionB}
The action of $\Gamma$  by automorphisms on $R$
extends in an action of $\Gamma$ by automorphisms on $B=R(\!(x;d)\!)$ if and only
if it is $d$-compatible. We have then:
\begin{equation}\label{actionsurB}
x^{-1}\cdot\gamma=p_{\gamma}x^{-1}+p_{\gamma}r_{\gamma} \ \ \text{
for any } \gamma\in\Gamma,
\end{equation} where $p\in Z^1(\Gamma,U(R))$ is the 1-cocycle associated to the $d$-compatibility and $r$ is an arbitrary
map $\Gamma\to R$ satisfying the identity:
\begin{equation}\label{relationsurr}
r_{\gamma\gamma'}=r_{\gamma'}+p_{\gamma'}^{-1}(r_{\gamma}\cdot\gamma')
\quad \text{ for all } \gamma,\gamma'\in\Gamma.\end{equation}
In particular, this action extends in
an action on $B$ if and only if it extends in an action on the
subring $A=R[x^{-1} ; -d]$.\end{theo}

\pf Suppose that $\Gamma$ acts by automorphisms on $B$ with $R \cdot \gamma=R$ for any $\gamma \in \Gamma$. 
We can apply point (i) of
proposition \ref{continuity} to write
$x^{-1}\cdot\gamma=g_{-1}x^{-1}+g_{0}+\sum_{j \geq 1}g_j x^j$, with $g_j\in R$
for any $j\geq -1$ and $g_{-1}\not=0$. Moreover $x^{-1}\in U(B)$
implies $(x^{-1}\cdot\gamma)\in U(B)$ which is equivalent to $g_{-1}\in U(R)$.
Applying $\gamma$ to \eqref{schur3}, we obtain:
$(x^{-1}\cdot\gamma)(f\cdot\gamma)-(f\cdot \gamma)(x^{-1}\cdot\gamma) =-d(f)\cdot\gamma$
for any $f\in R$. Since $f\cdot\gamma\in R$, we can develop this
identity:
$$[g_{-1}x^{-1}(f\cdot\gamma)-(f\cdot\gamma)g_{-1}x^{-1}]+
[g_0(f\cdot\gamma)-(f\cdot\gamma)g_0]+\sum_{j\geq
  1}[g_jx^j(f\cdot\gamma)-(f\cdot\gamma)g_jx^j]
=-d(f)\cdot\gamma.$$ The first term is:
$g_{-1}[x^{-1}(f\cdot\gamma)-(f\cdot\gamma)x^{-1}]=-g_{-1}d(f\cdot\gamma) \in
R$. The second term is zero by commutativity of $R$. The third term is of
valuation $\geq 1$. So we deduce
that:$$-g_{-1}d(f\cdot\gamma)=-d(f)\cdot\gamma\ \ \text{and} \ \
\sum_{j\geq 1}[g_jx^j(f\cdot\gamma)-(f\cdot\gamma)g_jx^j] =0 \ \ \text{for any } f\in R.$$

Denoting $p_{\gamma}=g_{-1}$, we have $p_{\gamma}\in U(R)$ and the
first equality means by definition that the action of $\Gamma$
is $d$-compatible. We claim that the second assertion implies $g_j=0$
for all $j\geq 1$. Suppose that there exists a minimal index
$m\geq 1$ such that $g_m\not=0$. Calculating with relation
\eqref{PDOxi}, we have $\sum_{j\geq
m}[g_jx^j(f\cdot\gamma)-(f\cdot\gamma)g_jx^j]=0$, then $mg_md(f\cdot\gamma)x^{m+1}+\cdots=0$ 
and therefore $d(f\cdot\gamma)=0$ for
any $f\in R$. It follows by $d$-compatibility of the action that $d=0$, hence
a contradiction. We
have proved that $x^{-1}\cdot \gamma=g_{-1}x^{-1}+g_0$ with
$p_{\gamma}=g_{-1}\in U(R)$ satisfying definition \ref{troisdef} (i). Then we
set $r_{\gamma}=(g_{-1})^{-1}g_0$. We have
$x^{-1}\cdot\gamma=p_{\gamma}x^{-1}+p_{\gamma}r_{\gamma}$. Relations
\eqref{cocycle}  for $p$ and  \eqref{relationsurr} for $r$ follow from relation
$(x^{-1}\cdot\gamma)\cdot\gamma'=x^{-1}\cdot\gamma\gamma'$.
\medskip

Conversely, let us assume that the action of $\Gamma$ on $R$ is
$d$-compatible. Denote by $p$ the 1-cocycle associated to the $d$-compatibility.
Let us choose a map $r:\Gamma\to R$ satisfying
\eqref{relationsurr}. For any $\gamma\in\Gamma$, we denote 
$q_{\gamma}=p_{\gamma}r_{\gamma}$ and calculate for all $f\in R$:
\medskip

\centerline{$(p_{\gamma}
x^{-1}+q_{\gamma})(f\cdot\gamma)-(f\cdot\gamma)(p_{\gamma}x^{-1}+q_{\gamma})
=p_{\gamma}(x^{-1}(f\cdot\gamma)-(f\cdot\gamma)x^{-1})=-p_{\gamma}d(f\cdot\gamma)$.}\medskip

Using \eqref{dcomp} we obtain 
$(p_{\gamma}x^{-1}+q_{\gamma})(f\cdot\gamma)-(f\cdot\gamma)(p_{\gamma}x^{-1}+q_{\gamma})=-d(f)\cdot\gamma$ for any $f\in R$.
Hence by \eqref{schur3} we obtain an action of $\Gamma$ by automorphisms on the polynomial algebra $A$ 
extending the initial action on $R$ by setting: $x^{-1}\cdot\gamma=p_{\gamma}x^{-1}+p_{\gamma}r_{\gamma}$.
The element $p_{\gamma}+q_{\gamma}x$ is invertible in $B_0$ and the element
$x^{-1}\cdot\gamma=p_{\gamma}x^{-1}+q_{\gamma}$ is invertible in $B$ because $p_{\gamma}\in U(R)$. 
Then we define  an action of $\Gamma$ by automorphisms on $B$ extending the action on $A$
by setting: $x\cdot\gamma=(x^{-1}\cdot\gamma)^{-1}=x(p_{\gamma}+q_{\gamma}x)^{-1}$.
The condition $(u\cdot\gamma)\cdot\gamma'=u\cdot\gamma\gamma'$ for any $\gamma,\gamma'\in\Gamma$ and any $u\in B$
follows from \eqref{cocycle} for $p$ and \eqref{relationsurr} for $r$.\epf

\begin{coro}\label{extensionBbis}
If the action of $\Gamma$  by automorphisms on $R$
is $d$-compatible, then it extends in an action of $\Gamma$ by automorphisms
on $B$ defined by:
\begin{equation}\label{actionsurBbis}
x^{-1}\cdot\gamma=p_{\gamma}x^{-1}, \ \text{ or equivalently }  \
x\cdot\gamma=xp_{\gamma}^{-1}=\sum_{j\geq
0}d^j(p_{\gamma}^{-1})x^{j+1} \ \text{ for any }
\gamma\in\Gamma,\end{equation}where $p\in Z^1(\Gamma,U(R))$ is 
the 1-cocycle associated to the $d$-compatibility.\end{coro}

{\it Proof.} We just apply theorem \ref{extensionB} for the
trivial map $r$ defined by $r_{\gamma}=0$ for any
$\gamma\in\Gamma$ which satisfies obviously
\eqref{relationsurr}.\qed

\begin{exas}\label{examplesr} We consider a $d$-compatible action of
$\Gamma$ on $R$. Let $p$ be the 1-cocycle associated  to the $d$-compatibility.
Relation \eqref{relationsurr} can be interpreted
as a 1-cocycle condition for the right action of $\Gamma$ on $R$ defined by
$\langle f\,|\,\gamma\rangle=p_{\gamma}^{-1}(f\cdot\gamma)$ for any
$\gamma\in\Gamma$ and $f\in R$. We denote by $Z^1_p(\Gamma,R)$ the
additive group of maps $r:\Gamma\to R$ satisfying
\eqref{relationsurr}. We consider here various examples for the
choice of $r\in Z^1_p(\Gamma,R)$.
\begin{itemize}
\item[1.] The case $r=0$ corresponds to the extension described in
corollary \ref{extensionBbis}. If $r$ is a coboundary (\ie there
exists $f\in R$ such that:
{$r_{\gamma}=p_{\gamma}^{-1}(f\cdot\gamma)-f$ \ for any $\gamma\in
\Gamma$}), then we can suppose up to a change of variables that
$r=0$, because the element $x'=(x^{-1}-f)^{-1}$ satisfies
$B=R(\!(x'\,;\,d)\!)=R(\!(x\,;\,d)\!)$ and
$x'^{-1}\cdot\gamma=p_{\gamma}x'^{-1}$ for any $\gamma\in\Gamma$.
\item[2.] A straightforward calculation proves that the map $r:\Gamma\to
R$ defined by: {$r_{\gamma}=-p_{\gamma}^{-1}d(p_{\gamma})$\ for
any $\gamma\in\Gamma$}  is an element of $Z^1_p(\Gamma,R)$. The
corresponding action of $\Gamma$ on $B$ is given by:
{$x^{-1}\cdot\gamma=p_{\gamma}x^{-1}-d(p_{\gamma})=x^{-1}p_{\gamma}$
\ for any $\gamma\in\Gamma$.}  \item[3.]  Since $Z^1_p(\Gamma,R)$
is a left $R^{\Gamma}$-module,  $\kappa\,r$ is an element
of $Z^1_p(\Gamma,R)$ for any $r\in Z^1_p(\Gamma,R)$ and any
$\kappa\in R^{\Gamma}$. The corresponding action of $\Gamma$ on
$B$ is given by: {$x^{-1}\cdot\gamma=p_{\gamma}x^{-1}+\kappa
\,p_{\gamma}r_{\gamma}$ for any $\gamma\in \Gamma$.} If we suppose
moreover that $\kappa\in U(R)$, then $x'=(\kappa^{-1}x^{-1})^{-1}$
satisfies $B=R(\!(x\,;\,d)\!)=R(\!(x'\,;\,\kappa^{-1}d)\!)$, and we obtain
$x'^{-1}\cdot\gamma=p_{\gamma}x'^{-1}+p_{\gamma}r_{\gamma}$ for any
$\gamma\in\Gamma$. Up to a change of variables, theses cases where $\kappa \in U(R)$ reduce to the case $\kappa=1$.
\end{itemize}
\end{exas}

In order to provide for $C$ an extension theorem similar to
\ref{extensionB}, we need the following technical lemma about square roots in $C$.

\begin{lemma}\label{squareroot} Let $q=\sum_{i\geq 2}g_iy^i$ be an element of 
$C$, with $g_i\in R$ for any $i\geq 2$, $g_2\not=0$. We suppose
that there exists $e\in U(R)$ such that $g_2=e^2$.\begin{itemize} \item[\rm{(i)}]
 There exists a unique element $z\in C$ of the form $z=ey+\sum_{i\geq
2}e_iy^i$ with $e_i\in R$ for any $i\geq 2$, such that $q=z^2$. \item[\rm{(ii)}] The only other series $z'$ satisfying $z'\,^2=q$ is $z'=-z$. \item[\rm{(iii)}] Moreover, if
$q\in B$, then $z\in By$.
\end{itemize}\end{lemma}

{\it Proof.} We compute inductively the coefficients $e_i$ of $z$
for $i\geq 2$ by identification in the equality $\sum_{i\geq 2}g_iy^i=(ey+\sum_{j\geq
2}e_iy^i)^2$. Using \eqref{QPDOyi} in the
development of the right hand side, we observe that
$z^2=e^2y^2+(2e_2e)y^3+(2e_3e+e_2^2+e\delta(e))y^4+\cdots=
e^2 y^2 + \sum_{i\geq 3}(2ee_{i-1}+h_{i-2})y^i$ where the rest $h_{i-2} \in R$ depends
only on previous elements $e,e_2,\ldots,e_{i-2}$ and their images
by $\delta$. Therefore $e_{i-1}=(2e)^{-1}(g_i-h_{i-2})$ and the proof
of (i) follows by induction.\medskip

Since $R$ is a domain, the only other element $e' \in R$ satisfying $e'\ ^2=g_2$ is $e'=-e$, then point \rm{(ii)} follows from point \rm{(i)}. \medskip

By point (iii) of \ref{BC}, we have $z=u+u'y$ with
$u,u'\in B$, $v_y(u)\geq 2$, $v_y(u')=0$. Suppose that $u\not=0$.
We denote $u=\sum_{i\geq m}f_{2i}y^{2i}$ and $u'=\sum_{i\geq
0}f'_{2i}y^{2i}$, with $m\geq 1$, $f_{2i},f'_{2i}\in R$, $f_{2m}\not=0$
and $f'_0=e$. Then $q=(u+u'y)^2=
(u^2+u'yu'y)+(uu'y+u'yu)$ and the assumption $q\in B$ implies
$uu'y+u'yu=0$ using point (iii) of proposition \ref{BC}. By \eqref{QPDOy1} and \eqref{QPDOy2} we have:\smallskip

\centerline{
$uu'y+u'yu=(f_{2m}y^{2m}+\cdots)(ey+\cdots)+(ey+\cdots)(f_{2m}y^{2m}+\cdots)=
2f_{2m}ey^{2m+1}+\cdots$}\smallskip

with $2f_{2m}e\not=0$. Hence a contradiction. \qed

\begin{theo}\label{extensionC}
The action of $\Gamma$  by automorphisms on $R$
extends in an action of $\Gamma$ by automorphisms on $C=R(\!(y ; \frac 12 d)\!)_2$ if and only
if it is quadratically $d$-compatible. In this case $B=R(\!(x ; d)\!)$ is stable under the extended action.\end{theo}

{\it Proof.} Suppose that the action extends in an action by automorphisms on
$C$. By point (ii) of proposition \ref{continuity}, we can
consider the map $\Gamma\to U(R)\,,\, \gamma\mapsto s_{\gamma}$
defined by $y\cdot\gamma=s_{\gamma}^{-1}y+\cdots$. Writing
$(y\cdot\gamma)\cdot\gamma'=(s_\gamma^{-1}\cdot\gamma')s_{\gamma'}^{-1}y+\cdots$, 
we observe with \eqref{cocycle} that $s\in Z^1(\Gamma,U(R))$. 
Moreover it follows from point (iii)
of proposition \ref{continuity} that the action restricts in an
action of $B$. Since
$x^{-1}\cdot\gamma=(y^{-1}\cdot\gamma)^2=(s_{\gamma}y^{-1}+\cdots)^2=s_{\gamma}^2x^{-1}+\cdots$,
we deduce from theorem \ref{extensionB} that $s_{\gamma}$
satisfies condition \eqref{quaddcomp}, therefore the action
is quadratically $d$-compatible.\medskip

Conversely, we assume now that the action of $\Gamma$ on $R$ is quadratically 
$d$-compatible. There exists some $s\in Z^1(\Gamma,U(R))$
satisfying \eqref{quaddcomp}. Then the map $p:\gamma\mapsto s_{\gamma}^2$ lies in $Z^1(\Gamma,U(R))$
and satisfies \eqref{dcomp}. Let $r$ be any map $\Gamma\to R$
satisfying \eqref{relationsurr}. By theorem \ref{extensionB}, we
define an action by automorphisms on $B$ by setting
$x^{-1}\cdot\gamma=p_{\gamma}x^{-1}+p_{\gamma}r_{\gamma}$ for any
$\gamma\in\Gamma$. The development of its
inverse in $B$ is of the form $x\cdot\gamma=s_{\gamma}^{-2}x+\cdots$. By lemma
\ref{squareroot}, there exists a unique element $\overline{y}_{\gamma}\in C$
such that $\overline{y}_{\gamma}^2=x\cdot\gamma$ and whose development is of the form $\overline{y}_{\gamma}=s_{\gamma}^{-1}y+\cdots$. In particular it follows from the
action of $\gamma$ on relation \eqref{QPDOy2} that:
\begin{equation}\label{temp}
\overline{y}_{\gamma}^2f=f\overline{y}_{\gamma}^2+\sum_{k \geq
2}2^{k-1}((\delta^{k-1}(f\cdot\gamma^{-1}))\cdot\gamma)\,\overline{y}_{\gamma}^{2k} \ \ \text{  for any } f\in R.
\end{equation}
We extend the action of $\gamma$ to $C$ by setting $y\cdot\gamma=\overline{y}_{\gamma}$.
% In order to extend the action of $\gamma$ to $C$ we set $y\cdot\gamma=\overline{y}_{\gamma}$
% and more generally:\begin{equation}\label{extactionC}q\cdot\gamma=\sum_{j\geq m}(f_j\cdot\gamma)\,\overline
% y^j \ \ \ \text{ for any } q=\sum_{j\geq m}f_jy^j\in C \text{ with} \ m\in\Z, \ f_j\in R.\end{equation} 
To prove that $q\mapsto q\cdot\gamma$ defines a automorphism
of $C$, it is sufficient by \eqref{QPDOy1} to prove that:
\begin{equation}\label{temp2}{ \overline{y}_{\gamma}f=f\overline{y}_{\gamma}+\sum_{k \geq
1}{\textstyle{\frac{(2k)!}{2^k(k!)^2}}}\,((\delta^{k}(f\cdot\gamma^{-1}))\cdot\gamma)\,
\overline{y}_{\gamma}^{2k+1}} \ \ \text{  for any } f\in
R.\end{equation} Since $\overline{y}_{\gamma}=s_{\gamma}^{-1}y+\cdots$ with $s_{\gamma}^{-1}\in U(R)$, any element of $C$ can be
written as a series in the variable $\overline{y}_{\gamma}$ with
coefficients in $R$. In particular, for any $f\in R$, we have
$\overline{y}_{\gamma}f=f\overline{y}_{\gamma}+\sum_{n\geq 1}\delta_n(f)\overline{y}_{\gamma}^{n+1}$, where $(\delta_n)_{n\geq 1}$ is a sequence of additive maps
$R\to R$. Hence $\overline{y}_{\gamma}^2f=f\overline{y}_{\gamma}^2+\sum_{n\geq
2}\Delta_n(f)\overline{y}_{\gamma}^{n+1}$ with notation
$\Delta_n=\sum_{j=0}^{n-1}\delta_j\delta_{n-j-1}$. The
identification with \eqref{temp} leads to $\Delta_n=0$ for even
$n$, and $\Delta_{n}=(2 \delta')^{k-1}$ for odd $n=2k-1$
where $\delta'$ is defined by $\delta'(f)=\delta(f\cdot\gamma^{-1})\cdot\gamma$ for any $f\in R$. It
follows by a straightforward induction that $\delta_n=0$ for any
odd $n$, then
$\Delta_{2k-1}=\sum_{i+j=k-1}\delta_{2j}\delta_{2i}$, and
finally $\delta_{2k}={\frac{(2k)!}{2^k(k!)^2}}\delta'^k$.
So relation \eqref{temp2} is satisfied. 
Because $(s_\gamma^{-1}\cdot\gamma')s_{\gamma'}^{-1}=s_{\gamma\gamma'}^{-1}$ 
for all $\gamma,\gamma'\in\Gamma$, we deduce from
$(x\cdot\gamma)\cdot\gamma'=x\cdot\gamma\gamma'$ with lemma \ref{squareroot} that 
$(y\cdot\gamma)\cdot\gamma'=y\cdot\gamma\gamma'$.
We conclude that it defines an action of $\Gamma$ by automorphisms on $C$.\qed

\begin{coro}\label{extensionCbis} 
We suppose that the action of $\Gamma$  by automorphisms on $R$
is quadratically $d$-compatible.
\begin{itemize}
\item[\rm{(i)}] Let $s\in Z^1(\Gamma,U(R))$ be a 1-cocycle associated to the quadratic $d$-compatibility. For $\gamma \in \Gamma$, let $\overline{y}_{\gamma}$ be the square root of $xs_{\gamma}^{-2}$  whose coefficient of minimal
valuation in its development as a series in the variable $y$ is $s_{\gamma}^{-1}$.  Then 
the action extends in an action of $\Gamma$ by automorphisms
on $C$ defined from
\begin{equation}\label{actionsurCbis}
y\cdot\gamma=\overline{y}_{\gamma}=s_{\gamma}^{-1}y+\cdots\ \text{ for any }
\gamma\in\Gamma\end{equation}
\item[\rm{(ii)}] The other
extensions of the action are given by
$y\cdot\gamma=\epsilon_{\gamma}\overline{y}_{\gamma}$ for any $\gamma\in
\Gamma$, where $\epsilon$ is a multiplicative map
$\Gamma\to\{-1,+1\}$.\end{itemize}\end{coro}

{\it Proof.} We apply the second part of the proof of theorem \ref{extensionC} to the
case where $r$ is defined by $r_{\gamma}=0$ for any
$\gamma\in\Gamma$.\qed

%%%%%%%%%%%%%%%%%%%%%%%%%%%%%%%%%%%%%%%%%%%%%%%%%%%%%%%%%%%%%%%%%%%%%%%%%%%%%%%%%%%%%%%%%%%%%%%%%%%%%%%%
%%%%%%%%%%%%%%%%%%%%%%%%%%%%%%%%%%%%%%%%%%%%%%%%%%%%%%%%%%%%%%%%%%%%%%%%%%%%%%%%%%%%%%%%%%%%%%%%%%%%%%%%

\subsection{Invariants in $B$ and $C$ for the extensions of actions on $R$}

We take all data and notations of proposition \ref{BC}.
For $\Gamma$ a group acting
by automorphisms on $B$ (respectively on $C$) stabilizing $R$,
we give sufficient conditions for the invariant ring $B^{\Gamma}$
(respectively $C^{\Gamma}$) to be described as a ring of pseudodifferential
(respectively quadratic pseudodifferential) operators with
coefficients in $R^{\Gamma}$. 

\begin{theo}\label{casimportantB}Let $\Gamma$ be a group acting by automorphisms on $B=R(\!(x;d)\!)$ stabilizing
$R$. We assume that there exists  in $B^{\Gamma}$ an element of the form 
$w=gx^{-1}+h$  with $h \in R$ and $g\in U(R)$. Then the derivation $D=g d$
 restricts to a derivation of $R^{\Gamma}$,
 and we have $A^{\Gamma} = R^{\Gamma}[w ; -D]$, 
$B_0^{\Gamma} = R^{\Gamma}[[w^{-1} ; D]]$ and $B^{\Gamma} = R^{\Gamma}(\!(w^{-1} ; D)\!)$.
\end{theo}

{\it Proof.} We consider the subring
$A=R[x^{-1} ; -d] \subset B$ and commutation law \eqref{schur3}.
We denote by $D$ the derivation of $R$ defined by $D(f) = g d(f) =
-(gx^{-1}f-fgx^{-1}) = -(wf-fw)$. Then $A=R[w ; -D]$. 
For any $f\in R^\Gamma$, $wf-fw\in R^\Gamma$ because $w\in B^{\Gamma}$.
Hence the restriction of $D$ to $R^\Gamma$ is a derivation of $R^\Gamma$.
Since $w$ is invariant, an element
of $A$ written as a polynomial in $w$ with coefficients in $R$ is
invariant if and only if any coefficient lies in $R^{\Gamma}$. We
conclude that $A^{\Gamma} = R^{\Gamma}[w ; -D]$.\smallskip

The element $u=w^{-1} \in B_0^{\Gamma}$ satisfies $v_x(u)=1$ and its dominant
coefficient $g^{-1}$ is invertible in $R$. 
We have $x=g(1+z)u$ for some $z\in B_0$ satisfying $v_x(z)\geq 1$.
Then $u$ is a uniformizer of the valuation $v_x$ in $B_0$, which means that any element of $B_0$ 
can be written as a power series in the variable $u$ with coefficients in $R$.
In particular, for $q=\sum_{i\geq 0}f_ix^i$ an element of $B_0^\Gamma$ with $f_i\in R$,
it follows from point (i) of proposition \ref{continuity} that $f_0\in R^\Gamma$ and
$q-f_0\in B_0^\Gamma$. We have $q-f_0=q'u$ with 
$q'=(\sum_{i \geq 0} f_{i+1}x^i)g(1+z) \in B_0$. Since
$q, f_0$ and $u$ belong to $B_0^{\Gamma}$ we deduce that $q' \in
B_0^{\Gamma}$. We have proved that for any $q \in B_0^{\Gamma}$,
there exist $f_0 \in R^{\Gamma}$ and $q' \in B_0^{\Gamma}$ such
that $q = f_0+q'u$. Applying this process for $q'$, there
exist $f'_0 \in R^{\Gamma}$ and $q'' \in B_0^{\Gamma}$ such that
$q=f_0+f'_0u+q''u^2$. It follows by induction that $q$ lies
in the left $R^\Gamma$-module $R^{\Gamma} \langle \! \langle u \rangle\!\rangle$
of power series in the variable $u$ with coefficients in $R^\Gamma$.
We conclude that $B_0^{\Gamma} \subset R^{\Gamma} \langle \! \langle u \rangle\!\rangle$.
The converse inclusion is clear, so $B_0^{\Gamma}=R^{\Gamma} \langle \! \langle u \rangle\!\rangle$. In particular, $R^{\Gamma} \langle \! \langle u \rangle\!\rangle$ is a subring of $B_0$. Hence for any $f\in R^\Gamma$, there exist
a sequence $(\delta_n(f))_{n\geq 0}$ of elements of $R^\Gamma$ such that
$uf=\sum_{n \geq 0} \delta_n(f) u^{n+1}$. The commutation
relation $wf-fw=-D(f)$ becomes $uf-fu=uD(f)u$. We compute:
$uf=fu+uD(f)u=fu+[D(f)u+uD^2(f)u]u=fu+D(f)u^2+[D^2(f)u+uD^3(f)u]u^2$,
and conclude by iteration that $\delta_n(f)=D^n(f)$ for any $n\geq 0$. In other words,
$B_0^{\Gamma} = R^{\Gamma}[[u ; D]]$.\smallskip

Let $R^{\Gamma}(\!(u ; D)\!)$ be the localized ring of
$B_0^{\Gamma}$ with respect of the powers of $u$.
Since $u$ is invertible in $B$, we have $R^{\Gamma}(\!(u ; D)\!)
\subset B$, and therefore $R^{\Gamma}(\!(u ; D)\!) \subset
B^{\Gamma}$. Conversely for any $f \in B^{\Gamma}$, there exists
an integer $n\geq 0$ such that $fu^n \in B_0$; since $fu^n \in
B_0^{\Gamma}=R^{\Gamma}[[u;D]]$, we deduce that $f \in
R^{\Gamma}(\!(u ; D)\!)$. Hence $B^{\Gamma}=R^{\Gamma}(\!(u ; D)\!)$ and
the proof is complete.\qed\bigskip

\begin{theo}\label{casimportantC}Let $\Gamma$ be a group acting by automorphisms on $C=R(\!(y;\delta)\!)_2$ stabilizing
$R$. Denote by $s$ the 1-cocycle in $Z^1(\Gamma,U(R))$ defined by:
 $y\cdot\gamma=s_{\gamma}^{-1}y+\cdots$ for any $\gamma\in \Gamma$.
We assume that there exists in $C^{\Gamma}$ an element of the
form $w=e^2y^{-2}+h$  with $h\in R$ and $e\in U(R)$ such that $e\cdot\gamma=s_{\gamma}^{-1}e$
for any $\gamma\in\Gamma$.
 Then the derivation $\Delta=e^2\delta$
 restricts to a derivation of $R^{\Gamma}$, and denoting by $v$ the square root
 of $w^{-1}$ whose development starts with $v=e^{-1}y+\cdots$, we
 have $C^{\Gamma}=R^{\Gamma}(\!(v\,;\,\Delta)\!)_2$.
\end{theo}

{\it Proof.}  Let us recall that the existence of the map $s$ follows from point (ii) of proposition \ref{continuity}. We know by point (iii) of proposition \ref{continuity}
that $B$ is stable under the action of $\Gamma$ on $C$.
We can apply theorem \ref{casimportantB} with $g=e^2$ to deduce that
$D=e^2d=2e^2\delta$
restricts to a derivation of $R^{\Gamma}$ and $B^{\Gamma}=R^{\Gamma}[[u\,;\,D]]$
where $u=w^{-1}\in B_0^{\Gamma}$. Since $u=e^{-2}y^2+\cdots$ lies in $B$, it follows from lemma
\ref{squareroot} that there exist in $C$ two elements $v=\sum_{j\geq 1}g_jy^j$ and
$v'=-v$ such that $v^2=v'^2=u$, with notations $g_j\in R$ for any $j\geq 1$, $g_j=0$ for any even
index $j$ and $g_1=e^{-1}$. For any $\gamma\in\Gamma$, we
have: $u=u\cdot\gamma=(v\cdot\gamma)^2$. Then $v\cdot\gamma=\pm v$. Since $v\cdot\gamma=(e^{-1}y+\cdots)\cdot\gamma=
s_{\gamma}e^{-1}(s_{\gamma}^{-1}y+\cdots)+\cdots=e^{-1}y+\cdots$, we are necessarily in the case $v\cdot\gamma=v$.
Therefore $v\in C^{\Gamma}$.\medskip

Let us denote $v=e^{-1}(1+z)y$ where $z=\sum_{j\geq 2}eg_jy^{j-1}\in C_0$ with $v_y(z)\geq 1$.
Similarly there exists $z'\in C_0$ with $v_y(z')\geq 1$ such that $y=e(1+z')v$.
Then $v$ is a uniformizer of the valuation $v_y$ in $C_0$, which means that any element of $C_0$ 
can be written as a power series in the variable $v$ with coefficients in $R$.
Using inductively  point (ii) of proposition \ref{continuity},
we can prove on the same way as in the proof of \ref{casimportantB} 
that $C_0^{\Gamma}=R^{\Gamma}\langle\!\langle v\rangle\!\rangle$. In particular, $R^{\Gamma}\langle\!\langle v\rangle\!\rangle$ is a subring of $C_0$. Hence for any $f\in R^\Gamma$, there exists
a sequence $(\delta_n(f))_{n\geq 0}$ of elements of $R^\Gamma$ such that
$vf=\sum_{n\geq 0}\delta_n(f)v^{n+1}$.
Comparing with relation $v^{-2}f-fv^{-2}=-D(f)$ which follows from theorem \ref{casimportantB},  direct calculations show that
$\delta_i=0$ for any odd index $i$, and $\delta_{2k}=\frac{(2k)!}{2^k(k!)^2}(\frac D2)^k$.
In other words $C^{\Gamma}=R^{\Gamma}(\!(v\,;\,\Delta)\!)_2$ with notation $\Delta=\frac12D$.\qed

\begin{coro}\label{schema} Under the assumptions of the previous theorem, and defining in
$R^{\Gamma}$ the derivations $\Delta=e^2\delta$ and $D=2\Delta=e^2d$, we have the following ring
embeddings:\medskip

\centerline{ \xymatrix{ A=R[x^{-1}\,;\,-d] \ \ar@{^{(}->}[r]
\ & \ B=R(\!(x\,;\,d)\!) \
\ar@{^{(}->}[r]_{x=y^2} \ & \ \ C=R(\!(y\,;\,\delta)\!)_2 \ \\
A^{\Gamma}=R^{\Gamma}[w\,;-D] \
\ar@{^{(}->}[r]\ar@{^{(}->}[u] \ & \
B^{\Gamma}=R^{\Gamma}(\!(w^{-1}\,;\,D)\!) \
\ar@{^{(}->}[r]_{w^{-1}=v^2}\ar@{^{(}->}[u] \
& \ C^{\Gamma}=R^{\Gamma}(\!(v\,;\,\Delta)\!)_2\ar@{^{(}->}[u] \\
}}
\end{coro}

{\it Proof.} It is a joint formulation of theorems \ref{casimportantB} and \ref{casimportantC}.\qed

%%%%%%%%%%%%%%%%%%%%%%%%%%%%%%%%%%%%%%%%%%%%%%%%%%%%%%%%%%%%%%%%%%%%%%%%%%%%%%%%%%%%%%%%%%%%%%%%%%%%%%%%
%%%%%%%%%%%%%%%%%%%%%%%%%%%%%%%%%%%%%%%%%%%%%%%%%%%%%%%%%%%%%%%%%%%%%%%%%%%%%%%%%%%%%%%%%%%%%%%%%%%%%%%%

\subsection{Weighted invariants in $R$ and equivariant splitting maps}\label{splittingproblem}

The data and notations are those of \ref{mainnota}. We introduce here for any quadratically $d$-compatible action of a group $\Gamma$ on $R$ a
natural link between the invariants in $C$ or $B$ and some weighted invariants in $R$.

\begin{defis}\label{weightinv} 
For any $1$-cocycle $s$ in $Z^1(\Gamma,U(R))$ and any integer $k$ we define:
\begin{equation}\label{weightkaction}
\trsf fk{\gamma}=s_{\gamma}^{-k}(f\cdot\gamma)
\quad \text{for any }f\in R, \ \gamma\in\Gamma.\end{equation} It follows from
relations \eqref{actionautom} and \eqref{cocycle} that:
\begin{equation}
\trsf{\trsf fk{\gamma}}k{\gamma'}=\trsf fk{\gamma\gamma'}\quad\text{for any } f\in R, \
\gamma,\gamma'\in\Gamma.\end{equation}
Hence relation \eqref{weightkaction} defines a right action of
$\Gamma$ on $R$ named the weight $k$ action associated to $s$. In
particular the weight zero action is just the original action of
$\Gamma$ by automorphisms on $R$:\begin{equation}\label{weight0action}
\trsf f0{\gamma}=f\cdot\gamma
\quad \text{for any }f\in R, \ \gamma\in\Gamma.\end{equation} We introduce the additive
subgroup of weight $k$ invariants:
\begin{equation}\label{weightkinv}
M_k=\{f\in R\,;\,\trsf fk\gamma=f \text{ for any
}\gamma\in\Gamma \}=\{f\in R\,;\,f\cdot\gamma=s_{\gamma}^kf \ \text{
for any }\gamma\in\Gamma \}.
\end{equation}
We have $M_0=R^{\Gamma}$ and $M_kM_{\ell}\subseteq M_{k+\ell}$ for
all $k,\ell\in\Z$.\end{defis}

\begin{nota}\label{Bgammar} 
We suppose that the action of $\Gamma$ on $R$ is quadratically $d$-compatible.
Let $s\in Z^1(\Gamma,U(R))$ be a 1-cocycle associated to the quadratic $d$-compatibility. 
We define $p=s^2\in Z^1(\Gamma,U(R))$ and we choose a map $r:\Gamma\to R$ satisfying
\eqref{relationsurr}. We consider the action of $\Gamma$ on $B$ and $C$
determined by theorems \ref{extensionB} and \ref{extensionC}. We
denote by $B^{\Gamma}$ and $C^{\Gamma}$ the corresponding
invariant rings. For any integer $k$, we introduce:
\begin{equation*}
C_k=\{q\in C\,;\,v_y(q)\geq k\},\hspace{1cm}  C_k^{\Gamma}=C_k\cap C^{\Gamma},
\end{equation*}
\begin{equation*}
B_k=\{q\in B\,;\,v_x(q)\geq k\},\hspace{1cm} B_k^{\Gamma}=B_k\cap B^{\Gamma}.
\end{equation*}
\end{nota}

\begin{prop}\label{pik} Let $k$ be an integer.
\begin{itemize}\item[\rm(i)] 
Let $\pi_k:C_k\to R$ be the canonical projection $\sum_{i\geq k}f_iy^i\mapsto f_k$.
The restriction of $\pi_k$ to $C_k^\Gamma$ defines an additive map $\pi_k:C_k^{\Gamma}\to M_k$.\item[\rm(ii)] 
Let $\widehat \pi_{k}:B_k\to R$ be the canonical projection $\sum_{i\geq k}f_ix^i\mapsto f_k$.
Then $\widehat \pi_k$ is the restriction of $\pi_{2k}$ to $B_{k}$, and its restriction to
$B_{k}^\Gamma$ defines an additive map $\widehat\pi_{k}:B_k^{\Gamma}\to
M_{2k}$.\end{itemize}\end{prop}

\pf For any $q=\sum_{i\geq k}f_iy^i\in C_k$ with $f_i\in
R$, $f_k\not=0$, and for any $\gamma\in\Gamma$, it follows
from theorem \ref{extensionC} that:
$q\cdot\gamma=(f_k\cdot\gamma)(s_{\gamma}^{-1}y+\cdots)^k+(f_{k+1}\cdot\gamma)(s_{\gamma}^{-1}y+\cdots)^{k+1}+\cdots$.
Then $q\cdot\gamma=(f_k\cdot\gamma)s_{\gamma}^{-k}y^k+\cdots$. Hence
$q\cdot\gamma=q$ implies $(f_k\cdot\gamma)s_{\gamma}^{-k}=f_k$, or equivalently
$f_k\in M_k$ by \eqref{weightkinv}. Point \rm(ii) follows from the second part of theorem \ref{extensionC}.\epf

\begin{problem}\label{splittingpb}Point (i) of this proposition leads to the natural question of finding
additive right splitting maps $\psi_k:M_k\to C_k^{\Gamma}$ such that
$\pi_k\circ \psi_k=\id_{M_k}$. 
Solutions arise by restriction if we can construct
$\psi_k:R\to C_k$ satisfying $\pi_k\circ \psi_k=\id_R$ and the equivariance
condition:\begin{equation}\label{equationsplitting}\psi_k(\trsf fk\gamma)
=\psi_k(f)\cdot\gamma\quad\text{ for any } \gamma\in \Gamma,f\in R.\end{equation}
In the particular case of even weights the corresponding question of finding splitting maps $\psi_{2k}:M_{2k}\to B_k^{\Gamma}$ was solved in \cite{CMZ} and \cite{UU} in various contexts involving modular forms.\end{problem}

\section{Application to algebraic modular forms}\label{secondsection}
\subsection{Homographic action on $R$}\label{21}

\begin{nota}\label{datamodular}
In order to apply the previous algebraic results in the
arithmetical context of modular forms, we specialize the
data and notations of \ref{mainnota}. From now $\Gamma$ is a
subgroup of ${\rm SL}(2, \C)$, and $R$ is a commutative
$\C$-algebra of functions in one variable $z$, which is a domain. We suppose in the following that:
\begin{itemize}
\item[(i)] $\Gamma$ acts on the right by homographic
automorphisms on $R$: \begin{equation}\label{defhomo}(f\cdot\gamma)(z)=
f\left(\frac{az+b}{cz+d}\right) \text{ for any }f \in R \text{ and }\gamma =
\left(\begin{smallmatrix} a & b \\ c & d \end{smallmatrix}\right) \in \Gamma,\end{equation} 
\item[(ii)] $z \mapsto cz+d \in U(R)$
for any $\gamma = \left(\begin{smallmatrix} a & b \\ c & d
\end{smallmatrix}\right) \in \Gamma$,  \item[(iii)] $R$ is stable under the standard
derivation $\partial_z$ with respect to $z$.
\end{itemize}\end{nota}

\begin{exas}\label{modularexamples}
A formal algebraic example of such a situation is the case where $R=\C(z)$, the field of complex rational
functions in one indeterminate. Numbertheoretical examples can arise from the following 
construction. Assume that $\Gamma$ is a subgroup of ${\rm SL}(2, \R)$ and denote
by $\mathcal{H} = \{z \in \C ; \Ima(z) >0\}$ the Poincar\'e upper
half-plane. Then $\mathcal{H}$ is stable by the homographic action
of $\Gamma$, and various subalgebras $R$ of holomorphic or meromorphic
functions on $\mathcal{H}$ satisfy the previous conditions. For
instance:
\begin{enumerate}
 \item $R_1 = \text{Hol}(\mathcal{H})$ the algebra of holomorphic
 functions on $\mathcal{H}$.  \item $R_2 = \text{Mer}(\mathcal{H})$ the
 field of meromorphic functions on $\mathcal{H}$.  \item For
 $\Gamma \subset {\rm SL}(2, \Z)$, $R_3$ is the field of $f \in \text{Mer}(\mathcal{H})$ such that, for any $p \in \mathbb{P}_1(\Q)$, there exists $V_p$ a $\mathcal{H}$-neighborhood
 of $p$ such that $f$ has neither zero nor pole in $V_p$. Here we consider
 the hyperbolic topology on $\overline{\mathcal{H}} = \mathcal{H} \cup \mathbb{P}_1(\Q)$,
 where a fundamental system of $\mathcal{H}-$neighborhoods of $p \in \mathbb{P}_1(\Q)$ is
 given by the upper half-planes $\{\Ima(z) > M\}$ for $p = \infty$ and by the open
 disks $D_r = \{z \in \mathcal{H} ; |z-(p+ir)|<r\}$ for $p \in \Q$. We can prove (see for instance \cite{FB} theorem 2.1.1)
that $R_3^{{\rm SL}(2, \Z)}=\C(\!(j)\!)$ where $j$ is the modular invariant. It follows in particular that the field $R_3$ satisfies 
$R_3^\Gamma\not=\C$ for any subgroup $\Gamma$ of ${\rm SL}(2, \Z)$. 
\end{enumerate} \end{exas}

\begin{nota}\label{datamodular2}
With data and notations \ref{datamodular}, we introduce moreover, according to proposition \ref{BC}, the noncommutative $\C$-algebras:
\begin{equation}
A=R[x^{-1}\,;\,\partial_z]\subset B=R(\!(x\,;\,d_z)\!) \subset C=R(\!(y\,;\,\delta_z)\!) \quad\text{with \ } 
x=y^2, \ d_z=-\partial_z, \ \delta_z=\frac12d_z.
\end{equation}
\end{nota}

\begin{lemma}\label{quadcomp}
 The map $s: \Gamma \to U(R), \gamma\mapsto s_{\gamma}$
 defined by:
 \begin{equation}\label{defdes}s_{\gamma}(z)=cz+d \quad\text{ for any }
 \gamma = \left(\begin{smallmatrix} a & b \\ c & d
\end{smallmatrix}\right) \in \Gamma\end{equation} is a 1-cocycle for the action of $\Gamma$.
\end{lemma}\pf A straightforward calculation using \eqref{defhomo} proves that \eqref{cocycle} is satisfied.\epf

\begin{defis}
Since $s\in Z^1(\Gamma,U(R))$, we can apply \ref{weightinv} to
introduce for any integer $k$ the weight $k$ action of
$\Gamma$ on $R$ defined by:
\begin{equation}\label{defactionmodk}
\trsf
fk{\gamma}(z)=\textstyle{(cz+d)^{-k}f\left(\dfrac{az+b}{cz+d}\right)}\qquad
\text{for any } f\in R \text{ and } \  \gamma=\left(\begin{smallmatrix} a
&b\\c&d\end{smallmatrix}\right)\in\Gamma, \end{equation} and the
$\C$-vector space of algebraic weight $k$ modular forms on $R$:
\begin{equation}\label{defformodk}
M_{k}=\{f\in R \,;\ \trsf fk{\gamma}=f \ \text{for any
}\gamma\in\Gamma\}.\end{equation}In particular for $k=0$
\begin{equation}\label{actionzero} \trsf f0\gamma=f\cdot\gamma\quad 
\text{and}\quad M_0=R^{\Gamma}.\end{equation}
\end{defis}

\begin{rem}\label{MkMl}
We use in the following the notations $$\mathcal{M}_j=\prod_{k \geq j} M_k,\,\,\, \mathcal{M}_j^{\ev}=\prod_{k \geq j} M_{2k} \text{ for any } j \in \Z \text{ \,\,and\,\, } \mathcal{M}_*=\bigcup_{j \in \Z}\mathcal{M}_j,\,\,\, \mathcal{M}_*^{\ev}=\bigcup_{j \in \Z}\mathcal{M}_j^{\ev}$$ 
with the convention $\mathcal{M}_{j_1} \subset \mathcal{M}_{j_2}$ for $j_1 \geq j_2$. Throughout the rest of the paper, we suppose that $M_k \cap M_\ell = \{0\}$ for all
integers $k \neq \ell$. In the
classical situations, it is sufficient for this additional
hypothesis to assume that $\Gamma$ contains at least one matrix
$\left(\begin{smallmatrix} a & b \\ c & d
\end{smallmatrix}\right)$ such that $(c, d) \notin \{0\} \times
\mathbb{U}_{\infty}$. Observe that this
excludes the unipotent cases where $\Gamma \subset \{\left(\begin{smallmatrix} 1 & a \\
0 & 1\end{smallmatrix}\right); a \in \R \}$. Then an element $\widetilde{f}$ of $\mathcal{M}_j$ can be denoted unambiguously by $\widetilde{f}=\sum_{k \geq j}f_k$, where $f_k \in M_k$.
\end{rem}

\begin{rem}\label{jacobi1}
For $k$ and $m$ nonnegative integers, denote by $\mathcal{J}_{k, m}$ the space of algebraic Jacobi forms on $R$ of weight $k$ and index $m$, defined as functions $\Phi: \mathcal{H} \times \mathbb{C} \to \mathbb{C}$ satisfying the Jacobi transformation equation: $$\Phi(\frac{az+b}{cz+d}, \frac{Z}{cz+d})=(cz+d)^ke^{\frac{2i\pi mcZ^2}{cz+d}}\Phi(z, Z) \hspace{0.5cm} \text{ for any } \gamma = \left(\begin{smallmatrix} a & b \\ c & d \end{smallmatrix}\right) \in \Gamma,$$
 and admitting around $Z=0$
  a Taylor expansion $\Phi(z, Z)=\sum_{\nu \geq 0} X_{\nu}(z)Z^{\nu}$ with $X_{\nu} \in R$ for any $\nu \geq 0$. Note that $\mathcal{J}_{k, m}$ is isomorphic to $\mathcal{J}_{k, 1}$ for any $m \geq 1$ via the map $Z \mapsto \sqrt{m}Z$. Then, for any $k \geq 0$ and $m \geq 1$, the vector space $\mathcal{J}_{k, m}$ is isomorphic to $\mathcal{M}_k$ ; an isomorphism is explicitly described page 34 of \cite{EZ} where the space $\mathcal{J}_{k, m}$ is denoted by $M_{k, m}$. 
\end{rem}

%%%%%%%%%%%%%%%%%%%%%%%%%%%%%%%%%%%%%%%%%%%%%%%%%%%%%%%%%%%%%%%%%%%%%%%%%%%%%%%%%%%%%%%%%%%%%%%%%%%%%%%%%%%%%%%%%
%%%%%%%%%%%%%%%%%%%%%%%%%%%%%%%%%%%%%%%%%%%%%%%%%%%%%%%%%%%%%%%%%%%%%%%%%%%%%%%%%%%%%%%%%%%%%%%%%%%%%%%%%%%%%%%%%
%%%%%%%%%%%%%%%%%%%%%%%%%%%%%%%%%%%%%%%%%%%%%%%%%%%%%%%%%%%%%%%%%%%%%%%%%%%%%%%%%%%%%%%%%%%%%%%%%%%%%%%%%%%%%%%%%

\subsection{Homographic action on $B$ and $C$}\label{22}

The data and notations are those of \ref{datamodular2}.

\begin{lemma}\label{compmodu} The homographic action of $\Gamma$ on $R$ is quadratically $\partial_z$-compatible, and an associated 1-cocycle is the map $s$ defined by \eqref{defdes}.\end{lemma}
\pf We have already observed in example 2 of \ref{exemplesactions} that, for any $\gamma = \left(\begin{smallmatrix} a & b \\
c & d \end{smallmatrix}\right) \in \Gamma$ and any $f\in R$, we have
$\partial_z(f\cdot\gamma)(z) =(cz+d)^{-2}\partial_z f(\frac{az+b}{cz+d})$.
In other words, $\partial_z f\cdot\gamma=s_{\gamma}^2\partial_z
(f\cdot{\gamma})$. From definition \ref{troisdef} (ii) and lemma \ref{quadcomp}, the lemma is proved.\epf
We can then apply the results of the first part of the paper:
for the canonical choice $r=0$ (see example 1 of \ref{examplesr}) we obtain the following results.
\begin{prop}\label{extensionBCkappa0}{ \ }\begin{itemize}\item[{\rm(i)}]The homographic action of $\Gamma$
on $R$ extends in an action by automorphisms on $B$ defined by:
\begin{equation}\label{actmodsurB}\trsf{x^{-1}}{}{\gamma} =
(cz+d)^2x^{-1}\quad \text{ for any } \gamma=\left(\begin{smallmatrix} a & b \\ c & d
\end{smallmatrix}\right) \in
\Gamma.\end{equation}
Consequently, for any $q  = \sum_{n>-\infty} f_n
x^n \in B$ with $f_n\in R$, we have:
\begin{equation}\label{actmodsurBbis}
\trsf q {}{\gamma} = \sum_{n>-\infty} (f_n\cdot{\gamma})
\trsf{x^{-1}}{}{\gamma}^{-n}\quad \text{ for any } \gamma \in
\Gamma.\end{equation}
\item[{\rm(ii)}]The homographic action of $\Gamma$
on $R$ extends in an action by automorphisms on $C$ defined by:
\begin{equation}\label{actmodsurC}\trsf{y}{}{\gamma} =
(cz+d)^{-1}y+\cdots\quad \text{ for any } \gamma=\left(\begin{smallmatrix} a & b \\ c & d
\end{smallmatrix}\right) \in\Gamma,\end{equation}
where this Laurent series is the square root in $C$ of 
$\trsf{x}{}{\gamma}=x(cz+d)^{-2}$ whose term of minimal valuation is $(cz+d)^{-1}y$. Consequently, for any $q  = \sum_{n>-\infty} f_n
y^n \in C$ with $f_n\in R$, we have $\trsf q {}{\gamma}= \sum_{n>-\infty} (f_n\cdot{\gamma})
\trsf{y}{}{\gamma}^{n}$ for any $\gamma \in
\Gamma$.
\item[{\rm(iii)}] The subalgebra $B$ of $C$ is stable under the action \textnormal{(ii)}, and the action \textnormal{(i)} is the restriction to $B$ of the action \textnormal{(ii)}.
\end{itemize}\end{prop} 
\pf We apply theorem \ref{extensionB} and theorem \ref{extensionC}.\epf

The two following theorems give explicit formulas describing the action of $\Gamma$ on $B$ and $C$ introduced in proposition \ref{extensionBCkappa0}.

\begin{theo}\label{extensionBC}
The extension to $B$ and $C$ of the homographic action of $\Gamma$ on $R$ given by proposition \ref{extensionBCkappa0} satisfy:
\begin{equation}\label{gammax}
\trsf{x}{}{\gamma} =  \sum_{n \geq 0} (n+1)!(cz+d)^{-2}\left(\frac c{cz+d}\right)^n \!\!\!x^{n+1}
\end{equation}
 \begin{equation}\label{gammay}
 \trsf{y}{}{\gamma} = \sum_{u \geq 0} \frac{(2u+1)!(2u)!}{16^{u}(u)!^3} (cz+d)^{-1} \left(\frac c{cz+d}\right)^{u} y^{2u+1}
 \end{equation}
for any $\gamma = \left(\begin{smallmatrix} a & b \\ c & d \end{smallmatrix} \right) \in \Gamma$.
\end{theo}

\pf
We have, from relations \eqref{PDOx1} and \eqref{actmodsurB},
$$\trsf{x}{}{\gamma} \!= \!\trsf{x^{-1}}{}{\gamma}^{-1} \!= \!x(cz+d)^{-2} = \sum_{u\geq 1}d_z^{u-1}((cz+d)^{-2})x^{u} = \sum_{u \geq 1} \frac{u!}{c^2}\left(\frac c{cz+d}\right)^{u+1} \!\!\!x^{u}$$
which proves formula \eqref{gammax}. We find here the action of $\Gamma$ described in \cite{CMZ} formula 1.7.

In order to prove \eqref{gammay}, we introduce the development $\trsf{y}{}{\gamma} \!= \!\sum_{i \geq 1} a_iy^i$ in $C$ with $a_i \in R$. Using the identity $\trsf{y}{}{\gamma}^2=\trsf{x}{}{\gamma}$ we prove by induction on $i$ that $a_i=0$ for even $i$ and that there exists a sequence $(\rho_j)_{j \geq 1}$ of complex numbers such that $\trsf{y}{}{\gamma} = \sum_{j \geq 1} \frac{\rho_j}c \left(\frac c{cz+d}\right)^{j} y^{2j-1}.$
Then we have
$$\trsf{y}{}{\gamma}^2 = \sum_{n, m \geq 1} \frac{\rho_n \rho_m}{c^2} \X^{n}\left[y^{2n-1}\X^m\right] y^{2m-1}.$$
From \eqref{QPDOyi}, and using the combinatorial identity
\begin{equation}\label{prodimpairs}
 \prod_{i=0}^{s-1}(2m+1+2i) = 2^{-s}\frac{(2m+2s)!m!}{(m+s)!(2m)!} \, \,\,\, \text{for any } \, m \geq 0 \,\text{ and }\, s \geq 1,
\end{equation}
we have $y^{2n-1}f =  \sum_{\ell \geq n-1} \frac{(2\ell)!(n-1)!}{\ell ! (2n-2)!}\frac{d_z^{\ell-n+1}(f)}{4^{\ell-n+1}(\ell-n+1)!}y^{2\ell+1}$, and an easy computation gives $d_z^t\left(\X^m\right)=\frac{(m+t-1)!}{(m-1)!}\X^{m+t}$. We deduce that:
%\sum_{k \geq 0} \frac{(2k+2n-2)!(n-1)!}{(k+n-1)!(2n-2)!}\frac{d_z^{k}(f)}{4^k k!}y^{2k+2n-1}
\begin{eqnarray*}
 \trsf{y}{}{\gamma}^2&=&\!\!\!\sum_{n, m \geq 1} \frac{\rho_n \rho_m}{c^2}\!\sum_{\ell \geq n-1} \frac{(2\ell)!(n-1)!(m+\ell-n)!}{4^{\ell-n+1}(\ell-n+1)!(2n-2)!\ell!(m-1)!}\X^{\ell+m+1}\!\!\!y^{2\ell+2m} \\
 &=& \sum_{u \geq 1}\left(\sum_{m=1}^{u}\sum_{i=m}^u \frac{\rho_{u+1-i}\rho_m}{4^{i-m}(i-m)!}\frac{(2u-2m)!(u-i)!(i-1)!}{(2u-2i)!(u-m)!(m-1)!}\right)\frac 1{c^2}\X^{u+1}\!\!\!y^{2u},
\end{eqnarray*}
with the change of variables $u=\ell+m$, and $i=u+1-n$. By identification in the equality $\trsf{y}{}{\gamma}^2 = \trsf{x}{}{\gamma}$, the coefficients $\rho_j$ satisfy, for any $u \geq 1$, the equation:
\begin{equation}\label{formulepourlesrho}
 \sum_{m=1}^{u}\sum_{i=m}^u \frac{\rho_{u+1-i}\rho_m}{4^{i-m}(i-m)!}\frac{(2u-2m)!(u-i)!(i-1)!}{(2u-2i)!(u-m)!(m-1)!} = u!
\end{equation}
The relation for $u=1$ gives $\rho_1^2=1$, then the choice of sign for $\rho_1$ determines inductively all the terms $\rho_u$, $u \geq 1$. Hence relation \eqref{formulepourlesrho} determines uniquely up to sign the sequence $(\rho_u)_{u \geq 1}$. Therefore the proof will be complete if we check that the coefficients $\rho_u =  \frac{(2u-1)!(2u-2)!}{16^{u-1}(u-1)!^3}$ satisfy relation \eqref{formulepourlesrho}. Then we have to prove that for any $u \geq 1$:
\begin{equation}\label{verifrho}
\sum_{m=1}^u \frac{(2u-2m)!(2m-1)!(2m-2)!}{4^{2u+m-2}(m-1)!^4(u-m)!} A(u,m) = u!
\end{equation}
with notations $A(u, m)=\sum_{i=m}^u G(u, m, i)$ and  $G(u, m, i) = \frac{4^i(i-1)!(2u-2i+1)!}{(i-m)!(u-i)!^2}$. Using Zeilberger's algorithm (see \cite{A=B}) we obtain the following relation which is easy to check directly: 
$$(4u+6)G(u, m, i)-(u+1-m)G(u+1, m, i)=H(u, m, i+1)-H(u, m, i)$$
where $H(u, m, i) = \frac{4^i(i-1)!(2u+3-2i)!}{(i-m-1)!(u+1-i)!^2}$. The summation on $i$ gives the relation $(4u+6)A(u,m)-(u+1-m)A(u+1, m)=0$, so $A(u,m)=A(m,m) \prod_{t=m}^{u-1}\frac{2(2t+3)}{t+1-m}$. Using $A(m,m)=4^m(m-1)!$ and \eqref{prodimpairs} we obtain $A(u,m) = 4^m\frac{(m-1)!(2u+1)!m!}{(u-m)!(2m+1)!u!}$.
So equation \eqref{verifrho} is equivalent to $T(u)=1$, where $T(u)=\sum_{m=1}^u K(u,m)$  with $K(u,m)=\frac{(2u+1)!}{16^{u-1}u!^2}\frac{(2m-2)!(2m-1)!m!(2u-2m)!}{(2m+1)!(m-1)!^3(u-m)!^2}$. Using again Zeilberger's algorithm we find the relation $$K(u+1,m)-K(u,m)=J(u,m)-J(u,m+1)$$  with $J(u,m)=\frac{(2u-2m+2)!(2u+1)!(2m-2)!}{16^u u u!(u+1)!(m-1)!(m-2)!( u+1-m)!^2}$. The summation on $m$ gives $T(u+1)-T(u)=0$, so $T(u)=T(1)=1$ for any $u \geq 1$, and the proof is complete.\epf

The next step is to describe the action of $\Gamma$ on any power of $y$.
We proceed depending on the parity and the sign of the exponent.
In particular relations \eqref{gammaxk} and \eqref{gammax-k} below are related to the action \eqref{gammax} on $B$ and already
appeared in previous studies about modular forms of even weight ([CMZ], [UU]).

\begin{lemma}
For any $k \in \N$, we have:
\begin{equation}\label{gammax-k}
  \trsf{y^{-2k}}{}{\gamma} = \trsf{x^{-k}}{}{\gamma} = \sum_{u=0}^{k-1} u! \binom{k}{u}\binom{k-1}u (cz+d)^{2k} \left(\frac c{cz+d}\right)^{u} x^{u-k},
 \end{equation}
 \begin{equation}\label{gammaxk}
  \trsf{y^{2k}}{}{\gamma} = \trsf{x^k}{}{\gamma} = \sum_{u \geq 0} u! \binom{k+u-1}{u}\binom{k+u}u (cz+d)^{-2k} \left(\frac c{cz+d}\right)^{u} x^{u+k},
  \end{equation}
 \begin{equation}\label{gammay2k+1}
  \trsf{y^{2k+1}}{}{\gamma} = \frac{(k+1)!k!}{(2k+2)!(2k)!}\sum_{u \geq 0} \frac{(2k+2u)!(2k+2u+2)!}{16^u u!(k+u)!(k+u+1)!}(cz+d)^{-2k-1}
  \left(\frac c{cz+d}\right)^{u} y^{2k+1+2u},
 \end{equation}
 \begin{multline}\label{gammay-2k+1}
\trsf{y^{-2k+1}}{}{\gamma} = \frac{(2k)!(2k-2)!}{k!(k-1)!}(cz+d)^{2k-1}
\left[\sum_{u=0}^{k-1}\frac{(k-u)!(k-1-u)!}{16^u u!(2k-2u)!(2k-2-2u)!}(\frac{c}{cz+d})^uy^{-2k+1+2u} \right. \\
- \left. \sum_{u \geq 0} \frac{(2u)!(2u+2)!}{u!(u+1)!}\frac
  1{16^{u+k}(u+k)!}(\frac{c}{cz+d})^{u+k}y^{2u+1}\right].
\end{multline}\end{lemma}
\pf
For the negative even powers of $y$ (the negative powers of $x$), we prove relation \eqref{gammax-k} inductively using the formulas 
$\trsf{x^{-1}}{}{\gamma} = (cz+d)^2x^{-1}$ and $\trsf{x^{-k-1}}{}{\gamma} = \trsf{x^{-1}}{}{\gamma}\trsf{x^{-k}}{}{\gamma}$ for $k>0$ with the relation $x^{-1}f = fx^{-1}+\partial_zf$. 

For the positive even powers of $y$ (the positive powers of $x$) denote, for any $k \geq 1$, $\trsf{y^{2k}}{}{\gamma}=\sum_{n \geq 0} a_k(n)(cz+d)^{-2k-n}c^ny^{2k+2n}$ with $a_k(n) \in \mathbb{C}$. Equations $\trsf{x^{-1}}{}{\gamma} = (cz+d)^{-2}x^{-1}$ and  $\trsf{y^{2k}}{}{\gamma} = \trsf{x^{-1}}{}{\gamma}\trsf{y^{2k+2}}{}{\gamma}$ give the 
relations $a_{k+1}(n)=(2k+n+1)a_{k+1}(n-1)+a_k(n)$ for $n \geq 1$. Then a double induction on $k$ and $n$ proves relation \eqref{gammaxk} using the expression of $a_1(n)$ for any $n \geq 0$ and $a_k(0)=1$ for any $k>0$.

We consider now relation $\eqref{gammay2k+1}$. Let $k \geq 0$; using
$\eqref{gammaxk}$,  $\eqref{gammay}$ and the relation
$\trsf{y^{2k+1}}{}{\gamma} = \trsf{x^k}{}{\gamma}\trsf{y}{}{\gamma}$ we obtain: $\trsf{y^{2k+1}}{}{\gamma} = \sum_{s \geq 0} \alpha_k(s) (cz+d)^{-2k-1}\left(\frac c{cz+d}\right)^s y^{2k+1+2s}$,
where 
$  \alpha_k(s)=\sum_{r=0}^s \frac{(2r+1)!(2r)!(s-k+r-1)!}{16^rr!^4(k-1)!k!}\beta_{k, s}(r)$, with
$  \beta_{k, s}(r) = \sum_{j=0}^{s-r} \frac{(k+s-r-j)!(r+j)!}{(s-r-j)!j!}$.
By Zeilberger's algorithm we obtain $(k+s+2)\beta_{k, s}(r)-(s-r+1)\beta_{k, s+1}(r)=0$, which proves that $\beta_{k, s}(r) = \frac{k!r!(k+s+1)!}{(s-r)!(k+r+1)!}$; 
more precisely: $(k+s+2)b(s, j)-(s-r+1)b(s+1, j) = G(s, j+1)-G(s,j)$ where $b(s,j) = \frac{(k+s-r-j)!(r+j)!}{(s-r-j)!j!}$ and 
$G(s,j)=\frac{(k+s-r-j+1)!(r+j)!}{(j-1)!(s-r-j+1)!}$.
Applying again Zeilberger's algorithm to $C(s, r)=\frac{(2r+1)!(2r)!(k+s-r-1)!}{16^rr!^3(s-r)!(k+r+1)!}$,  we obtain: 
$$(2s+2k+3)(2s+2k+1)C(s, r)-4(s+1)(k+s+2)C(s+1, r)=H(s, r+1)-H(s,r),$$
with $H(s,r)=\frac{4(k+s-r)!(2s+1)!(2r)!}{16^r(s-r+1)!(k+r)!(r-1)!r!^2}$. 
It follows that $\alpha_k(s)=\frac{(k+1)!k!(2k+2s)!(2k+2s+2)!}{(2k+2)!(2k)!16^ss!(k+s)!(k+s+1)!}$, 
which proves formula $\eqref{gammay2k+1}$.

Finally, for formula \eqref{gammay-2k+1}, denote
$\trsf{y^{-2k+1}}{}{\gamma} = \sum_{j \geq 0} \alpha_{-k}(j)(cz+d)^{2k-1}(\frac{c}{cz+d})^j y^{-2k+1+2j}$ with $\alpha_{-k}(j) \in \mathbb{C}$. The equality 
$\trsf{y^{-2k-1}}{}{\gamma} = \trsf{x^{-1}}{}{\gamma}\trsf{y^{-2k+1}}{}{\gamma}$ gives for any nonnegative $j$ the relation
\begin{equation}\label{recuneg}
 \alpha_{-(k+1)}(j)=\alpha_{-k}(j)+(2k-j)\alpha_{-k}(j-1),
\end{equation}
with the convention $\alpha_{-k}(-1)=0$. We proceed then inductively on $k$: the expression of $\alpha_0(j)$ is given by equation $\eqref{gammay}$. 
We have $\alpha_{-k}(0)=1$ for any $k \geq 0$, and by \eqref{recuneg} $\alpha_{-1}(j)= -4j^2 \frac{(2j-1)!(2j-2)!}{16^j j!^3}$ for $j \geq 1$ proves formula \eqref{gammay-2k+1} for $k=1$ and gives the base case.

For the inductive step, we prove that $\alpha_{-(k+1)}(j)$ is equal to the corresponding term of formula  \eqref{gammay-2k+1} separating the three cases $j \leq k-1$, $j=k$ and $j \geq k+1$. We check by direct computations using induction hypothesis for $k$ that the right hand side of formula \eqref{recuneg} corresponds to the expected expression in formula \eqref{gammay-2k+1} for $\alpha_{-(k+1)}(j)$.
\epf

We are now in position to give an unified formula for the action of $\Gamma$ on any power of $y$. 

\begin{theo}\label{extensionyk}
 For any $k \in \Z$, we have: 
\begin{equation}\label{gammayk}
 \trsf{y^k}{}{\gamma} = \sum_{u \geq 0} \omega_k(u) (cz+d)^{-k} \left(\frac c{cz+d}\right)^{u} y^{2u+k}
 \end{equation}
with notation $\omega_k(0)=1$ \ and \ $  \omega_k(u) = \frac 1{4^uu!} \prod_{i=0}^{u-1}(k+2i)(k+2i+2)$ for any $u \geq 1$.
\end{theo}

\pf Using formula \eqref{prodimpairs}, we check that formula \eqref{gammayk} corresponds in each case ($k$ even or odd, positive or negative) to formulas \eqref{gammaxk}-\eqref{gammay-2k+1}. 
\epf

%%%%%%%%%%%%%%%%%%%%%%%%%%%%%%%%%%%%%%%%%%%%%%%%%%%%%%%%%%%%%%%%%%%%%%%%%%%%%%%%%%%%%%%%%%%%%%%%%%%%%%%%%%%%%%%%%
%%%%%%%%%%%%%%%%%%%%%%%%%%%%%%%%%%%%%%%%%%%%%%%%%%%%%%%%%%%%%%%%%%%%%%%%%%%%%%%%%%%%%%%%%%%%%%%%%%%%%%%%%%%%%%%%%
%%%%%%%%%%%%%%%%%%%%%%%%%%%%%%%%%%%%%%%%%%%%%%%%%%%%%%%%%%%%%%%%%%%%%%%%%%%%%%%%%%%%%%%%%%%%%%%%%%%%%%%%%%%%%%%%%

\subsection{Equivariant splitting maps and noncommutative product on modular forms}\label{23}
The data and notations are those of previous paragraphs \ref{21} and \ref{22}. In order to simplify the notations, we will denote $f^{(n)}=\partial^n_z f$ for any $f \in R$ and $n \geq 0$.
According to problem \ref{splittingpb} and in connection with the construction already known in the even case, we seek to construct for any $m\in\Z$ a linear morphism $\psi_m:R\to C$
satisfying the conditions:
\begin{itemize}
 \item[(C1)] $\psi_{m}(\trsf{f}{m}{\gamma})=\trsf{\psi_{m}(f)}{}{\gamma}$ for any $f\in R,\gamma\in\Gamma$;
 \item[(C2)] there exists some complex sequence $(\alpha_m(n))_{n \geq 0}$ such that $\psi_{m}(f) =  \sum_{n \geq 0} \alpha_m(n)f^{(n)}y^{m+2n}$ for any $f\in R$;
 \item[(C3)] $\alpha_m(0)=1$.
\end{itemize}
We prove in the following that such a map $\psi_m$ exists and is unique for $m\geq 0$ (proposition \ref{relevementpositif}), exists but is not unique
if $m$ is negative and even (remark \ref{defPsi2Barre}), and doesn't exist if $m$ is negative and odd (proposition \ref{casimpairnegatif}).

\begin{prop}\label{relevementpositif}
For any nonnegative integer $m$, the linear map $\psi_m:R\to C_0$ defined by:
\begin{equation}\label{Tpositif}
\psi_m(f)=\sum_{n \geq 0}\frac{(-1)^n}{4^nn!} \left(\prod_{i=0}^{n-1}\frac{(m+2i)(m+2i+2)}{(m+i)}\right) f^{(n)}y^{m+2n}\quad\text{for any } f\in R
\end{equation}
is the unique map satisfying the three conditions \textnormal{(C1), (C2)} and \textnormal{(C3)}.
\end{prop}

\pf
The aim is to prove for any nonnegative integer $k$ the following formulas: 
\begin{equation}\label{Tpositifpair}
\psi_{2k}(f)  =  \sum_{n\geq 0}
 \frac{(2k-1)! \,(n+k-1)!}{(n+2k-1)! \,(k-1)!} \binom{-k-1}{n} f^{(n)} x^{k+n} 
\end{equation}
\begin{equation}\label{Tpositifimpair}
\psi_{2k+1}(f) = \frac{k!^2}{(2k+1)!}\sum_{n \geq 0} \frac{(-1)^n}{16^n} \frac{(2k+1+2n)!(2k+2n)!}{n!(2k+n)!(k+n)!^2}f^{(n)} y^{2k+1+2n}
\end{equation}
For even weights, formula $\eqref{Tpositifpair}$ and the unicity follow from propositions 2 and 6 of \cite{CMZ} up to a multiplicative constant in order to assure condition (C3). We prove now formula \eqref{Tpositifimpair} and show that this is the unique map satisfying the conditions (C1), (C2) and (C3) for odd positive weights.

We use at first the following equality for $m \geq 0$ and $n \geq 0$ integers:
\begin{equation}\label{deriveeactionmodulaire}
 \left(\trsf fn{\gamma}\right)^{(m)}(z) = \sum_{r=0}^m \frac{(-1)^{m-r}m!(m+n-1)!}{r!(m-r)!(n-1+r)!(cz+d)^{n+2r}}\X^{m-r}\!\!\!\!(f^{(r)}\cdot \gamma)(z)
\end{equation}
that can be easily checked by induction on $m$. Let $(\alpha_{2k+1}(n))_{n \geq 0}$ be a sequence of complex numbers, and 
let $\phi_{2k+1}(f) =   \sum_{n \geq 0} \alpha_{2k+1}(n) f^{(n)}y^{2k+1+2n}$. On the one hand, we have 
$\trsf{\phi_{2k+1}(f)}{}{\gamma} = \sum_{n \geq 0} \alpha_{2k+1}(n) (f^{(n)}\cdot \gamma) \trsf{y^{2k+1+2n}}{}{\gamma}$.
Using \eqref{gammay2k+1} we obtain (with $u=r+n$):
\begin{multline*}
 \trsf{\phi_{2k+1}(f)}{}{\gamma}(z) = \sum_{u \geq 0}\sum_{n=0}^u \frac{\alpha_{2k+1}(n)(k+n+1)!(k+n)!(2k+2u)!(2k+2u+2)!}{(2k+2n)!(2k+2n+2)!(k+u)!(k+u+1)!16^{u-n}} \\
 \frac{(cz+d)^{-2k-1-2n}}{(u-n)!}{\X}^{u-n}(f^{(n)}\cdot \gamma)(z)y^{2k+1+2u}.
\end{multline*}

On the other hand, using \eqref{deriveeactionmodulaire} we have:
 \begin{multline*}
\phi_{2k+1}(f|_{2k+1}\gamma)(z) = \sum_{u \geq 0} \alpha_{2k+1}(u)
(f|_{2k+1}\gamma)^{(u)}(z) y^{2k+1+2u} 
\\= \sum_{u \geq 0}\sum_{n=0}^u \alpha_{2k+1}(u)
\frac{u!(u+2k)!(-1)^{u-n}}{n!(n+2k)!} \frac{(cz+d)^{-2k-1-2n}}{(u-n)!}{\X}^{u-n}(f^{(n)}\cdot \gamma)(z)y^{2k+1+2u}.
\end{multline*}

Then we deduce that condition (C1) holds for $\phi_{2k+1}$ if and only if the sequence $(\alpha_{2k+1}(n))_{n \geq 0}$ satisfies, for any $u \geq 0$ and any $0 \leq n \leq u$, the equality: 
\begin{equation}\label{alphaku}
\alpha_{2k+1}(u)\frac{u!(u+2k)!(-1)^{u-n}}{n!(n+2k)!} =
\alpha_{2k+1}(n)\frac{(k+n+1)!(k+n)!(2k+2u)!(2k+2u+2)!}{(2k+2n)!(2k+2n+2)!(k+u)!(k+u+1)!16^{u-n}}.
\end{equation}

This equation is true for $n=u$ for any $u \geq 0$.
Let $u \geq 1$, and consider $n=u-1$. Then relation \eqref{alphaku} is equivalent to
\begin{eqnarray*}
 -u(u+2k)\alpha_{2k+1}(u) &=& \frac{\alpha_{2k+1}(u-1)}{16}
\frac{(k+u)!(k+u-1)!(2k+2u)!(2k+2u+2)!}{(k+u)!(k+u+1)!(2k+2u)!(2k+2u-2)!}
\\
% &=& \frac{\alpha_{2k+1}(u-1)}{16}
% \frac{(2k+2u+2)(2k+2u+1)(2k+2u)(2k+2u-1)}{(k+u+1)(k+u)} \\
&=& \alpha_{2k+1}(u-1) (k+u+\tfrac 12)(k+u-\tfrac 12).
\end{eqnarray*}

Assuming that $\alpha_{2k+1}(0)=1$ to satisfy condition (C3), we obtain as a necessary condition that, for any $u \geq 1$, $\frac{\alpha_{2k+1}(u)}{\alpha_{2k+1}(u-1)}=-\frac{(u+k+\frac 12)(u+k-\frac 12)}{u(u+2k)}$. Using formula \eqref{prodimpairs}, it follows that:
\begin{equation*}
\alpha_{2k+1}(u) = (-1)^u \prod_{i=1}^u
\frac{(k+1/2+i)(k-1/2+i)}{i(2k+i)} =\frac{(-1)^u}{16^u} \frac{k!^2}{(2k+1)!}
\frac{(2k+2u+1)!(2k+2u)!}{u!(u+2k)!(u+k)!^2}. 
% \frac{(-1)^u}{4^u}
% \frac{\prod_{i=0}^{u-1}(2k+3+2i)\prod_{i=0}^{u-1}(2k+1+2i)}{u!\frac{(u+2k)!}{(2k)!}}
% \\
% &=& \frac{(-1)^u}{16^u} \frac{(2k)!}{u!(2k+u)!}
% \frac{(2k+2u+1)!k!}{(k+u)!(2k+1)!} \frac{(2k+2u)!k!}{(k+u)!(2k)!} \\
\end{equation*} 

It proves consequently the unicity of the sequence $(\alpha_{2k+1}(u))_{u\geq 0}$. 

Let now $u \geq 1$ and $0 \leq n \leq u$, and let $\alpha_{2k+1}(u)=\frac{(-1)^u}{16^u} \frac{k!^2}{(2k+1)!}
\frac{(2k+2u+1)!(2k+2u)!}{u!(u+2k)!(u+k)!^2}$. We have on the one hand:
\begin{eqnarray*}
\alpha_{2k+1}(u)\frac{u!(u+2k)!(-1)^{u-n}}{n!(u-n)!(n+2k)!} 
&=&
% \frac{(-1)^n}{16^u}\frac{k!^2}{(2k+1)!}
% \frac{u!(u+2k)!(2k+2u+1)!(2k+2u)!}{n!(u-n)!(n+2k)!u!(u+2k)!(u+k)!^2}
% \\
% &=& 
\frac{(-1)^n}{16^u}\frac{k!^2}{(2k+1)!}
\frac{(2k+2u+1)!(2k+2u)!}{n!(u-n)!(n+2k)!(u+k)!^2},
\end{eqnarray*}

and on the other hand:
\begin{multline*}
\alpha_{2k+1}(n)\frac{(k+n+1)!(k+n)!(2k+2u)!(2k+2u+2)!}{(2k+2n)!(2k+2n+2)!(k+u)!(k+u+1)!16^{u-n}(u-n)!}
\\
% =
% \frac{(-1)^n}{16^u}\frac{k!^2}{(2k+1)!}
% \frac{(2k+2u)!(2k+2n+1)!(n+k+1)!(2k+2u+2)!}{(2k+2n+2)!(n+k)!(k+u+1)!n!(u-n)!(n+2k)!(u+k)!} \\
= \frac{(-1)^n}{16^u}\frac{k!^2}{(2k+1)!}\frac{(2k+2u+1)!(2k+2u)!}{n!(u-n)!(n+2k)!(u+k)!^2}
\end{multline*}

which proves that the sequence $(\alpha_{2k+1}(u))_{u\geq 0}$ satisfies equation \eqref{alphaku} for any $u \geq 0$ and $0 \leq n \leq u$, and gives the existence of the lifting $\psi_{2k+1}$.
\epf

According to the general principle of quantization by deformation, we can now transfer to modular forms the noncommutative product on operators in $C_0^{\Gamma}$. The case of even weights was considered in \cite{CMZ} and \cite{UU}.
\begin{theo}\label{*positif} We use notations \ref{datamodular2} and \ref{MkMl} and we denote $\mathcal M_0=\prod_{j \geq 0} M_j$ and $\mathcal M^{\ev}_{0}=\prod_{j \geq 0} M_{2j}$.
We consider the algebras $B_0=R[[x\,;\,d_z]] \subset C_0=R[[y\,;\,\delta_z]]_2$.
\begin{enumerate}
 \item[{\rm (i)}] The map $\Psi: \mathcal M_0 \to C_0^{\Gamma}$ defined by  $\Psi(\widetilde f)=\sum_{m \geq 0} \psi_{m}(f_m)$ for any $\widetilde f=\sum_{m\geq 0}f_m$ 
 is a vector space isomorphism, and $\mathcal M_0$ is an associative $\C$-algebra for the noncommutative  product:
 \begin{equation}\label{*+}\widetilde f \star \widetilde g = {\Psi}^{-1}(\Psi(\widetilde f)\cdot\Psi(\widetilde g)) \text{\quad for all }\ \widetilde f,\widetilde g\in\mathcal M_0.\end{equation}
 \item[{\rm (ii)}] In particular for any modular forms $f,g$ of respective nonnegative weights $k$ and $\ell$ (even or odd), the product $f \star g$ in $\mathcal M_0$ is:
  \begin{equation}\label{*+mod}f \star g = {\Psi}^{-1}(\psi_k(f)\cdot \psi_ \ell(g))=\sum_{n \geq 0} \alpha_n(k,  \ell)[f, g]_n \in \mathcal M_0,\end{equation}
where the coefficients $\alpha_n(k,  \ell)$ are rational constants depending only on $k, \ell$ and $n$, and
$[f, g]_n=\sum_{j=0}^n (-1)^j
\binom{k+n-1}{n-j} \binom{\ell+n-1}{j} f^{(j)} g^{(n-j)} \in
M_{k+\ell+2n}$ is the $n$-th Rankin-Cohen bracket of $f$ and $g$.
 \item[{\rm (iii)}] The restriction of $\Psi$ to the subspace $\mathcal M^{\ev}_0$ determines 
 a vector space isomorphism $\Psi_2$ between  $\mathcal M^{\ev}_0$ and $B_0^\Gamma$.\end{enumerate}
\end{theo}
\pf It is clear that $\Psi$ is linear and injective. Let $q=  \sum_{m=m_0}^{+\infty} h_m y^m$ an element of valuation $m_0 \geq 0$ in $C_0^{\Gamma}$. Then $h_{m_0} \in M_{m_0}$, and by condition (C3) the element $q-\psi_{m_0}(h_{m_0})$ lies in $C_0^{\Gamma}$ and its valuation is greater than $m_0$. It follows by induction that $q \in \Psi(\mathcal M_{m_0})$. Then $\Psi$ is a vector space isomorphism, and point {\rm (i)} is obtained by transfer of structures. Point {\rm (ii)} is a consequence of condition (C2) satisfied by $\psi_m$, and of the property of Rankin-Cohen brackets to be (up to constant) the unique operator $\sum_{k=0}^n a_n f^{(k)}g^{(n-k)}$ with $a_n \in \C$ mapping $M_k \times M_l$ to $M_{k+l+2n}$ (see for instance the end of paragraph 1 of \cite{Z2}).  Point (iii) is clear by condition (C2).
\epf

\begin{coro}\label{jacobi2} The invariant subspace $C_k^{\Gamma}$ is for any $k \geq 0$ isomorphic to the space $\mathcal{J}_{k, m}$ of algebraic Jacobi forms on $R$ of weight $k$ and positive index $m$. In particular relation \eqref{*+} defines a structure of noncommutative algebra on $\mathcal{J}_{k, m}$.
\end{coro}
\pf
Follows from theorem \ref{*positif} by remark \ref{jacobi1}.
\epf

\begin{rem}\label{isoeven}
The isomorphisms ${\Psi}$ and ${\Psi}^{-1}$ are explicit ; we can separate even and odd weights cases. In the even case, relations between the coefficients
of an element $\widetilde f=\sum_{n \geq 1} f_{2n}$ of $\mathcal M^{\ev}_0$ and its image ${\Psi}_2(\widetilde f) =
\sum_{m \geq 1} h_m x^m$ in $B_0^{\Gamma}$ are (as proved in \cite{CMZ}):
\begin{eqnarray}h_m &=& \sum_{r=0}^{m-1}\frac{(m-1)!(2m-2r-1)!}{(m-r-1)!(2m-r-1)!}
\binom{-m+r-1}{r} f_{2m-2r}^{(r)},\label{formulehmkappa}\\
f_{2n} &=& \frac{(n-1)!}{(2n-2)!} \sum_{r=0}^{n-1}
\frac{(2n-2-r)!}{(n-1-r)!}\binom{n}{r}h_{n-r}^{(r)}, \label{formulef2n}\end{eqnarray}
and the coefficients $\alpha_n(2k,2\ell)$ are computed for instance in \cite{CMZ}, \cite{UU}, \cite{Yao}.

In the odd case, technical calculations give the following relations between the coefficients of $\widetilde f=\sum_{n \geq 0} f_{2n+1}$ of $ \prod_{n \geq 0} M_{2n+1}$ and its image $\Psi(\widetilde{f})=\sum_{m \geq 0} h_{2m+1}y^{2m+1}$:
\begin{eqnarray}h_{2m+1} &=& \frac{(2m+1)!(2m)!}{m!^2}\sum_{r=0}^{m}\frac{(-1)^r(m-r)!^2}{16^rr!(2m-2r+1)!(2m-r)!}f_{2m-2r+1}^{(r)}, \label{formulehmodd}\\
f_{2n+1} &=& \frac{(2n)!(2n+1)!}{(2n-1)!n!^2} \sum_{r=0}^{n}
\frac{(2n-1-r)!(n-r)!^2}{16^r(2n-2r)!r!(2n-2r+1)!}h_{2n-2r+1}^{(r)}. \label{formulef2n+1}
\end{eqnarray}
\end{rem}

\begin{rems}\label{defPsi2Barre}
\begin{enumerate}
 \item[{\rm (i)}] The coefficients $\alpha_m(n)= \frac{(-1)^n}{4^nn!} \prod_{i=0}^{n-1}\frac{(m+2i)(m+2i+2)}{(m+i)}$ appearing in proposition \ref{relevementpositif}
are also well defined for $m$ even and non-positive. For $n \geq \frac{-m}2$, $\alpha_m(n)=0$ since $\alpha_{-2k}(n)=\frac 1{4^n n!}\prod_{i=0}^{n-1}\frac{(2k-2i)(2k-2-2i)}{(2k-i)}$ vanishes for $n \geq k$ because of the term $i=k-1$.
Therefore relation \eqref{Tpositifpair} can be completed for $k>0$ by:
\begin{equation}\label{Tnegatifpair}
 \psi_{-2k}(f) =  \sum_{n=0}^{k}
\frac{k! \,(2k-n)!}{(k-n)! \,(2k)!}\binom{k-1}{n}f^{(n)}x^{-k+n}.\end{equation}
The map $\psi_{-2k}: R \to B$ satisfies the three conditions (C1), (C2) and (C3), see for instance \cite{CMZ}. Then we can build a vector space isomorphism $\overline{\Psi}_2:\mathcal M_{*}^{\ev}\to B^\Gamma$ whose restriction to $\mathcal M^{\ev}_0$ is the isomorphism $\Psi_2$ of point (iii) of theorem \ref{*positif}. 
\item[{\rm (ii)}] Observe however that the lifting $\psi_{-2k}$ is not the only map satisfying the three conditions (C1), (C2) and (C3). Indeed, using the well-known
Bol's identity:
\begin{equation}\label{Bolidentity}
 \trsf f{2-h}{\gamma}^{(h-1)} = \trsf {f^{(h-1)}}h{\gamma} \quad
\text{ for any } \gamma \in {\rm SL}(2, \C), \ f\in R, \ h>0,
\end{equation}
the operator $\phi_c: f
\mapsto \psi_{-2k}(f)+c\psi_{2k+2}(f^{(2k+1)})$ satisfies for any $c \in \C$ conditions (C2) and (C3), as well as the equivariancy condition (C1):
\begin{eqnarray*} \trsf{\phi_c(f)}{}{\gamma}_{} &
= & \trsf{\psi_{-2k}(f)}{}{\gamma}_{} +
c \trsf{\psi_{2k+2}(f^{(2k+1)})}{}{\gamma}_{} \\
&=& \psi_{-2k}\trsf{f}{-2k}{\gamma} + c
\psi_{2k+2}\trsf{f^{(2k+1)}}{2k+2}{\gamma} \\
&=& \psi_{-2k}\trsf{f}{-2k}{\gamma} +
c\psi_{2k+2} \left(\trsf f{-2k}{\gamma}^{(2k+1)}\right)
\\&=& \phi_c\trsf{f}{-2k}{\gamma}.
\end{eqnarray*}
The lifting map $\psi_{-2k}$ is canonical in the sense that it's the only one which is polynomial, that is of the form $  \psi_{-2k}(f)=\sum_{n =
0}^k \alpha_n f^{(n)}x^{n-k} \in A$ for any $f\in R$.
\end{enumerate}
\end{rems}

The coefficients $\alpha_m(n)$ introduced in remark \ref{defPsi2Barre} (i) are not defined for $m=-2k+1$ with $k>0$ (for $n \geq 2k$ the denominator vanishes without vanishing of the numerator). In this negative and odd case, we have
the following result.

\begin{prop}\label{casimpairnegatif} Let $k$ be a positive integer. There is no map $\psi_{-2k+1}:R\to C$ satisfying the three conditions \textnormal{(C1), (C2)} and \textnormal{(C3)}.
More precisely, the only map $\psi_{-2k+1}:R\to C$ satisfying conditions \textnormal{(C1)} and \textnormal{(C2)} is defined (up to constant) by $\psi_{-2k+1}(f)=\psi_{2k+1}(f^{(2k)})$
 for any $f\in R$, and it doesn't satisfy condition \textnormal{(C3)}.\end{prop}

\pf Relation \eqref{deriveeactionmodulaire} can be extended to the cases where $m \geq 0$ and $n \in \Z$ by: 
\begin{equation}\label{deriveenegative}
 (f|_n\gamma)^{(m)}(z) = \sum_{r=0}^m \frac{m!}{r!} \binom{m+n-1}{m-r}
\frac{(-c)^{m-r}}{(cz+d)^{n+m+r}}f^{(r)}(\gamma z)
\end{equation}
which gives Bol's identity in the particular case $m=1-n$. Let $k>0$, suppose that there exists a linear morphism $\psi_{-2k+1}: R \to C$ satisfying (C1) and (C2). Using \eqref{gammay-2k+1} and \eqref{deriveenegative}, we can show as in the proof of proposition \ref{relevementpositif} that the coefficients of the series $\psi_{-2k+1}(f)=\sum_{n \geq 0} \alpha_{-2k+1}(n) f^{(n)}y^{-2k+1+2n}$ must satisfy the equality $ \alpha_{-2k+1}(n) \frac{n!}{r!} \binom{n-2k}{n-r}(-1)^{n-r} =$
\begin{equation*}
\begin{cases} 
    \frac{\alpha_{-2k+1}(r)}{16^{n-r}(n-r)!}
    \frac{(2k-2r)!(2k-2r-2)!(k-n)!(k-n-1)!}{(k-r)!(k-r-1)!(2k-2n)!(2k-2n-2)!}  & \mbox{ if } r \leq n \leq k-1, \\ \noalign{\vspace{3pt}} -\frac{\alpha_{-2k+1}(r)}{16^{n-r}(n-r)!} \frac{(2k-2r)!(2k-2r-2)!(2n-2k)!(2n-2k+2)!}{(k-r)!(k-r-1)!(n-k)!(n-k+1)!} & \mbox{ if } r \leq k-1 < n, \\    \noalign{\vspace{3pt}}
\frac{\alpha_{-2k+1}(r)}{16^{n-r}(n-r)!}\frac{(r-k+1)!(r-k)!(2n-2k)!(2n-2k+2)!}{(2r-2k+2)!(2r-2k)!(n-k)!(n-k+1)!}
 & \mbox{ if } k \leq r \leq n.
\end{cases}
\end{equation*}
For $n=2k$ and $r<2k$ we deduce that $\alpha_{-2k+1}(r)=0$ (separating the cases $r>k$ and $r \leq k$). Then condition (C3) can not be satisfied. Moreover, if we fix $\alpha_{-2k+1}(2k)=1$, it is easy to show that we have $\alpha_{-2k+1}(n)=\alpha_{2k+1}(n-2k)$ for any $n \geq 2k$, which proves the second assertion of the proposition. \epf

%%%%%%%%%%%%%%%%%%%%%%%%%%%%%%%%%%%%%%%%%%%%%%%%%%%%%%%%%%%%%%%%%%%%%%%%%%%%%%%%%%%%%%%%%%%%%%%%%%%%%%%%%%%%%%%%%
%%%%%%%%%%%%%%%%%%%%%%%%%%%%%%%%%%%%%%%%%%%%%%%%%%%%%%%%%%%%%%%%%%%%%%%%%%%%%%%%%%%%%%%%%%%%%%%%%%%%%%%%%%%%%%%%%
%%%%%%%%%%%%%%%%%%%%%%%%%%%%%%%%%%%%%%%%%%%%%%%%%%%%%%%%%%%%%%%%%%%%%%%%%%%%%%%%%%%%%%%%%%%%%%%%%%%%%%%%%%%%%%%%%

\subsection{Invariants for the homographic action on $B$}\label{24}
The data and notations are those of previous paragraphs \ref{21}, \ref{22} and \ref{23}. 

\begin{theo}\label{invBmodu} We assume that there exists a weight 2 modular form $\chi \in R$ which is
invertible in $R$. Then $B^\Gamma=R^\Gamma(\!(u\,;\,D)\!)$ and $B_0^{\Gamma}=R^\Gamma[[u\,;\,D]]$ for $u=x\chi\in B^\Gamma$ and $D=-\chi^{-1}\partial_z$.
\end{theo}

\pf It is clear that $u=x\chi$ lies in $B^\Gamma$, since its inverse
$w=\chi^{-1}x^{-1}$ is invariant by \eqref{actmodsurB}. Then we apply theorem \ref{casimportantB} with $g=\chi^{-1}$, $h=0$
and $d=-\partial_z$.\epf

\begin{rem}It follows from this theorem that the invariant algebra $B^{\Gamma}$ is noncommutative 
if and only if the restriction of $D$ to $R^\Gamma$ is not trivial, or equivalently
$R^\Gamma\not=\C$. Since the product of any modular form of even weight $2k$ by $\chi^{-k}$
lies in $R^\Gamma$, this condition is also equivalent to $M_{2k}\not=\C\chi^k$.
In the degenerate case where $R^\Gamma=\C$,  the product $\star$ defined by \eqref{*+mod}  is commutative
and most of the following results become much more simple. We refer to \ref{modularexamples} to give examples of fields $R$ such that $R^{\Gamma}$ is different from $\C$, and containing an invertible weight 2 modular form.
\end{rem}

\begin{problem} The explicit description of $B^{\Gamma}$ given by theorem \ref{invBmodu}
and the isomorphism $\overline{\Psi}_2:\mathcal M_*^{\ev}\to B^{\Gamma}$ introduced in \ref{defPsi2Barre} lead naturally to the question of finding relations
between the terms of a sequence of modular forms in $\mathcal M_*^{\ev}$ and the coefficients in $R^{\Gamma}$ 
of its image in $B^{\Gamma}$ by $\overline{\Psi}_2$. We answer this problem  for positive weights,
that is for the isomorphism $\Psi_2:\mathcal M_0^{\ev}\to B_0^{\Gamma}$.
Our first goal is to calculate, for any invariant pseudodifferential operator $q \in B_0^\Gamma$,
the modular forms $f_{2m}\in M_{2m}$ appearing in ${\Psi}_2^{-1}(q )=\sum_{m\geq 0}f_{2m}$. For this we introduce the following notation 
for the powers of the generator $u$:
\begin{equation}\label{defgk}
\text{for any } k \in \N, \  {\Psi}_2^{-1} (u^k) = \sum_{n \geq k} g_{k, 2n},\ \text{ with }
g_{k, 2n} \in M_{2n}. 
\end{equation}
\end{problem}

\begin{prop} The assumptions are those of theorem \ref{invBmodu}. Let $q \in B_0^{\Gamma}$ and $(a_k)_{k \geq 0}$ be the sequence of elements of $R^{\Gamma}$ such that $q =\sum_{k \geq 0} a_k u^k$. Then the modular forms $f_{2m} \in M_{2m}$ appearing as the terms of ${\Psi}_2^{-1}(q) = \sum_{m\geq 0} f_{2m}$
are given by:
$$f_{2m} = (-1)^m \frac{(m-1)!}{(2m-2)!}\sum_{k=0}^m\sum_{n=k}^{m} (-1)^{n}\frac{(2n-1)!(m-n)!}{(n-1)!(m+n-1)}
\binom{m}{m-n} [a_k, g_{k, 2n}]_{m-n}.$$
\end{prop}
\pf
The expression of any element of $B_0^{\Gamma}$ as $q=\sum_{k \geq 0} a_ku^k$ with $a_k \in R^{\Gamma}$ follows directly from theorem \ref{invBmodu}. By theorem \ref{*positif} we have:
\begin{equation*} 
{\Psi}_2^{-1}(q ) = \sum_{k}{\Psi}_2^{-1}(a_ku^k)  = \sum_{k} a_k \star {\Psi}_2^{-1} (u^k) =  \sum_{k, n} a_k \star g_{k, 2n} = \sum_{k, n, r} \alpha_{r}(0, 2n) [a_k, g_{k, 2n}]_r.
\end{equation*}
Then $f_{2m} = \sum_{k=0}^m\sum_{n=k}^m \alpha_{m-n}(0, 2n) [a_k, g_{k, 2n}]_{m-n}$ for any $m \geq 0$. \smallskip

By formula \eqref{Tpositifpair}, we obtain  $a_k\cdot \psi_{2n}(g_{k, 2n}) = \sum_{t \geq n} h_{t} x^{t}$, with
$$h_{t} = \frac{(2n-1)!(t-1)!}{(t+n-1)!(n-1)!}\binom{-n-1}{t-n} a_k \left(g_{k, 2n}\right)^{(t-n)}.$$
We deduce with formula $\eqref{formulef2n}$ that, for any $m \geq n \geq k$: 
$$\alpha_{m-n}(0, 2n)[a_k, g_{k, 2n}]_{m-n}  =  \frac{(m-1)!}{(2m-2)!}
\sum_{r=0}^{m-1} \frac{(2m-2-r)!}{(m-1-r)!} \binom{m}r h_{m-r}^{(r)}.$$
The left hand side is equal to \begin{equation*}\alpha_{m-n}(0, 2n) \sum_{i=0}^{m-n} (-1)^i 
\binom{m-n-1}{m-n-i}\binom{m+n-1}i a_k^{(i)} g_{k, 2n}^{(m-n-i)}, \end{equation*}
and expanding the right hand side gives the formula
$$\frac{(m-1)!(2n-1)!}{(2m-2)!(n-1)!}\sum_{i=0}^{m-n} 
\frac{a_k^{(i)}g_{k, 2n}^{(m-n-i)}}{i!}
\left[\sum_{r=i}^{m-n} \frac{(2m-2-r)!r!}{(m-r+n-1)!(r-i)!}\binom{m}{r}\binom{-n-1}{m-r-n}\right].$$
Identifying the terms for $i=m-n$ of each side gives the identity
$$\alpha_{m-n}(0,2n) = (-1)^{m-n} \frac{(m-1)!(2n-1)!(m-n)!}{(2m-2)!(n-1)!(m+n-1)}\binom{m}{m-n}$$
which proves the proposition.\epf

The next step is to obtain an explicit expression of the modular forms $g_{k, 2n}$ as functions of $\chi$. 

\begin{prop}\label{calculcalcul}
 For any nonnegative integers $k$ and $i$, we have:
 \begin{equation}\label{gk}
g_{k, 2k+2i} = (-1)^i \frac{(k+i)!(k+i-1)!}{(2k+2i-2)!k!}\hspace{-8pt}\sum_{\begin{smallmatrix}
(t_1, \ldots, t_k) \in (\Z_{\geq 0})^k \\ t_1+\ldots +
t_k=i \end{smallmatrix}} 
\hspace{-16pt}\gamma_i(t_1, \ldots, t_k) \frac{\chi^{(t_1)}}{(t_1+1)!}\cdots \frac{\chi^{(t_k)}}{(t_k+1)!}\end{equation}
where the coefficients $\gamma_i(t_1, \ldots, t_k)$ are defined for 
$(t_1, \ldots, t_k) \in (\Z_{\geq 0})^k$ such that $\sum_{j=1}^k t_j=i$ by:
\begin{equation}\label{alpha2}
 \gamma_i(t_1, \ldots, t_k) = \sum_{r=0}^i (-1)^r \frac{(2k+2i-2-r)!}{(k+i-1-r)!}
\sum_{\begin{smallmatrix}b_1+\ldots+b_k=r
    \\ 0 \leq b_1 \leq t_1 \\ \ldots \\ 0 \leq b_k \leq
    t_k\end{smallmatrix}}\!\!\!\!\binom{t_1+1}{b_1}\ldots
\binom{t_k+1}{b_k}\end{equation}and satisfy:
\begin{equation}\label{alpha1}
\frac{\gamma_i(t_1, \ldots, t_k)}{(k+i-1)!}= \hspace{-8pt}\sum_{a_1+\ldots+a_k \leq k+i-1} 
\binom{k+i-2}{a_1+\ldots +a_k-1} \prod_{j=1}^k [(-1)^{t_j} \delta(a_j=0) + \delta(a_j=t_j+1)]
\end{equation}with notation $\delta$ for the Kronecker symbol.
% $\delta(P)=1$ if $P$ is true, $\delta(P)=0$ otherwise.
\end{prop}

\pf Firstly let us observe that the sum in the right hand side of \eqref{alpha1} contains
$2^k-2$ terms corresponding to all $k$-tuples $(a_1, a_2, \ldots, a_k)$ 
with $a_j \in \{0, t_j+1\}$ for any $0 \leq j \leq k$, 
except the two $k$-tuples $(0, 0, \ldots, 0)$ and $(t_1+1, t_2+1, \ldots, t_k+1)$.\medskip

The proofs of identities \eqref{gk}, \eqref{alpha2} and \eqref{alpha1} are straightforward but quite technical, 
and we give only the main steps omitting intentionally some details.
Applying relation \eqref{PDOxi} to $u=x\chi$ we prove by induction on $k$ that 
$$u^k = \sum_{n \geq 0} \left(\frac{(-1)^n(k+n)!}{k!} \sum_{s_1+\ldots +
s_k=n} 
\frac{\chi^{(s_1)}}{(s_1+1)!} \ldots \frac{\chi^{(s_k)}}{(s_k+1)!}\right)x^{k+n}.$$
Then we apply \eqref{formulef2n} to obtain:
$$g_{k, 2n} = \frac{(-1)^{n-k}n!(n-1)!}{(2n-2)!k!}\!\!\!\!\!\!\!\! \sum_{t_1+\ldots+t_k=n-k}\left(\sum_{r=0}^{n-k}
\frac{(-1)^r(2n-2-r)!}{(n-1-r)!}\beta_r(t_1, \ldots, t_k)\right)
\frac{\chi^{(t_1)}}{(t_1+1)!} \ldots \frac{\chi^{(t_k)}}{(t_k+1)!}$$
where, for any $0 \leq r \leq n-k$,
$$\beta_r(t_1, \ldots, t_k) = \!\!\!\!\!\!\sum_{\begin{smallmatrix}
b_1+\ldots+b_k=r \\ b_1 \neq t_1+1, \ldots, b_k \neq t_k+1 \end{smallmatrix}}\!\!\!\!\!\!\binom{t_1+1}{b_1}\ldots \binom{t_k+1}{b_k},$$
which proves \eqref{gk} and \eqref{alpha2} with $\gamma_i(t_1, \ldots, t_k) = \sum_{r=0}^i(-1)^r \frac{(2k+2i-2-r)!}{(k+i-1-r)!}\beta_r(t_1, \ldots, t_k)$. Observe that $\beta_r(t_1, \ldots, t_k)$
is the $X^r$-coefficient in the polynomial $A(X) = \prod_{j=1}^k [(X+1)^{t_j+1}-X^{t_j+1}]$. We have $\gamma_i(t_1, \ldots, t_k) = Q^{(k+i-1)}(1)$ with notation 
$$Q(X)=\sum_{r=0}^i (-1)^r \beta_r(t_1, \ldots, t_k)X^{2k+2i-2-r} = X^{2i+2k-2}A(-\frac 1X).$$
We deduce that:
$$\gamma_i(t_1, \ldots, t_k) = 
\left.{\left[X^{k+i-2}\prod_{j=1}^k
\left((X-1)^{t_j+1}-(-1)^{t_j+1}\right)\right]^{(k+i-1)}}\right| _{X=1 },$$
which leads to \eqref{alpha1} using the relations: 
$$\left.\left[X^{k+i-2}\right]^{(b)}\right|_{X=1} = \left\{\begin{array}{cc}
\tfrac{(k+i-2)!}{(k+i-2-b)!}=b!\binom{k+i-2}{b} & \mbox{ if } b \leq k+i-2, \\ 0
& \mbox{ if }b=k+i-1,\end{array}\right.$$
and
$$\left.\left[(X-1)^{t_j+1}-(-1)^{t_j+1}\right]^{(a_j)}\right|_{X=1} = a_j!\left[(-1)^{t_j}\delta_{a_j=0} + \delta_{a_j = t_j+1}\right].$$
\epf

\begin{coro}\label{noformodeven}
 We have $g_{k, 2k+2i} = 0$ for any nonnegative integer $k$ and any nonnegative odd integer $i$.
\end{coro}

\pf We use the expression of the $g_{k, 2k+2i}$ given in the previous proposition.
Let $i$ an odd integer, and $(t_1, \ldots, t_k) \in (\Z_{\geq 0})^k$ such that 
$t_1+\ldots+t_k=i$. In the right-hand side of \eqref{alpha1}, 
we can group into pairs the $k$-tuples $(a_1, \ldots, a_k)$ and $(t_1+1-a_1,\ldots,t_k+1-a_k)$. The associated summands 
are opposite numbers because $i$ is odd. Hence $\gamma_i(t_1, \ldots, t_k)=0$.\epf

\begin{rem} The result obtained in corollary \ref{noformodeven} as a consequence of proposition \ref{calculcalcul} can also
be proved by theoretical arguments using algebraic properties of the
product $\star$ defined in \eqref{*+mod} and some combinatorial identity
conjectured in relation  (3.4) of \cite{CMZ} (see \cite{Yao} for a proof). We proceed by induction on $k \geq 2$. For $k=2$, applying \eqref{*+mod} 
with notations \eqref{defgk} we have:
$$\Psi_2^{-1}(u^2)=\chi \star \chi =\sum_{n \geq 0} \alpha_n(2, 2) [\chi, \chi]_n,$$
which proves by $(-1)^i$-symmetry of 
the Rankin-Cohen brackets that $g_{2, 4+2i} = \alpha_i(2,2)[\chi,\chi]_i = 0$ for any odd $i$. 
Suppose now that $\Psi_2^{-1}(u^k) = \sum_{i \in 2\Z_{\geq 0}} g_{k, 2k+2i}$ with $g_{k, 2k+2i} \in M_{2k+2i}$.
Then we compute: 
\begin{equation*} \Psi_2^{-1}(u^{k+1}) = 
\Psi_2^{-1}(u^{k}\cdot u) = \sum_{n \geq k,  \ 2|n-k} g_{k, 2n} \star \chi 
= \sum_{n \geq k,  \ 2|n-k}  \ \sum_{m \geq 0} \alpha_m(2n, 2)[g_{k, 2n}, \chi]_m
\end{equation*}
and similarly:
\begin{equation*} \Psi_2^{-1}(u^{k+1}) =
\Psi_2^{-1}(u \cdot u^{k}) = \sum_{n \geq k, \ 2|n-k} \chi \star g_{k, 2n} 
= \sum_{n \geq k,  \ 2|n-k} \ \sum_{m \geq 0} \alpha_m(2n, 2)[\chi, g_{k, 2n}]_m  
\end{equation*}
% &=& \sum_{n \geq k,  \ 2|n-k} \ \sum_{m \geq 0} \alpha_m(2n, 2) (-1)^m [g_{k, 2n}, \chi]_m.

Hence we deduce that, for each $r \geq k+1$:
$$\sum_{n+m=r}\alpha_m(2n, 2)[g_{k, 2n}, \chi]_m = \sum_{n+m=r}\alpha_m(2, 2n)(-1)^m[g_{k, 2n}, \chi]_m = g_{k+1, 2r}.$$
It follows from formula (3.4) of \cite{CMZ} that $\alpha_m(2, 2n) = \alpha_m(2n, 2)$. 
Using this identity and the induction hypothesis $g_{k, 2n}=0$ if $n-k$ is odd, we conclude that:
$$g_{k+1, 2r} = \sum_{\substack{m=0\\  m \text{ odd }}}^{r-k}\alpha_m(2r-2m, 2)[g_{k, 2r-2m}, \chi]_m 
= -\sum_{\substack{m=0\\  m \text{ odd }}}^{r-k} \alpha_m(2, 2r-2m)[g_{k, 2r-2m}, \chi]_m = 0.$$\end{rem}

\begin{rem} For even nonnegative integers $i$, the modular forms $g_{k, 2k+2i}$ defined by \eqref{defgk} are
described in proposition \ref{calculcalcul} as multidifferential polynomials of the fundamental weight two modular form $\chi$. Then we can necessarily express them in terms of Rankin-Cohen brackets. For small values of $i$, we compute:
$$g_{k, 2k} = \chi^k,$$
$$g_{k, 2k+4} = \frac{(k+2)!}{72(2k+1)(k-2)!}\,\chi^{k-2}[\chi, \chi]_2,$$
$$g_{k, 2k+8} = 2\frac{(2k+1)!(k+4)!(k+5)!}{(2k+6)!(k-1)!(k-2)!}\,\chi^{k-4} \left(\frac{[\chi, \chi^3]_4}{16920} + 
\frac{47k^2-187k+282}{121824k}[\chi, \chi]_2^2\right).$$
\end{rem}

\subsection{Invariants for the homographic action on $C$}\label{25}
The data and notations are those of previous paragraphs \ref{21}, \ref{22} and \ref{23}. 
A natural question is whether there exists for the action of $\Gamma$ on $C$ studied in paragraph \ref{22}
a description of the invariant subalgebra $C^\Gamma$ similar to the one obtained in theorem \ref{invBmodu}
 for $B^\Gamma$. 
%  The first step (point \rm(ii) of the following theorem) is the straightforward application of theorem \ref{casimportantC}
% to deduce some analogue of theorem \ref{invBmodu}. 
% A main difficulty for explicit calculations
% is then the complicated form of the square root $v$. 
% Another approach could be to introduce the element $z=\psi_1(\xi)^{-1}\in C^\Gamma$ as a generator of
% $C^\Gamma$ , but in this case the commutation law in the Laurent power series $C^\Gamma$
% in the indeterminate $z$ with coefficients in $R^\Gamma$ is not of the type \eqref{PDOx1} or \eqref{QPDOy1}, but twisted by a general higher derivation
% (in the sense of \cite{FD2}). That is why the detailed results on the invariants are limited to the framework of pseudodifferential operators.\medskip

\begin{theo}\label{invCmodu} We suppose that there exists a weight 1 modular form
$\xi$ which is invertible in $R$. Then we have $C^\Gamma=R^\Gamma(\!(v\,;\,\Delta)\!)_2$ for $\Delta=-\frac12\xi^{-2}\partial_z$ and $v=\xi y+\cdots$ a square root of $y^2\xi^2$.
\end{theo}

\pf It is the direct application in the modular situation of theorem \ref{casimportantC}.
\epf

\begin{rem}
 The previous theorem describes the elements of \(C^{\Gamma}=R^\Gamma(\!(v\,;\,\Delta)\!)_2\) as series 
with coefficients in \(R^\Gamma\) where the uniformizer \(v\) is choosen as
a square root of \(u=x\xi^2=\psi_2(\xi^2)\). In particular \(C^\Gamma\)
contains the ring of differential operators \(A^\Gamma=R^\Gamma[u^{-1}\,;\,-2\Delta]\).
The main obstacle to explicit calculations in this case is the complicated shape of the
square root \(v\).

Another more natural idea would be to consider the uniformizer \(z=\psi_1(\xi)\).
By reasoning as in theorem \ref{casimportantC}, we can prove then that \(C^{\Gamma}=
R^{\Gamma}(\!(z\,;\,S)\!)\) where the product $zf$ with 
\(f\in R^\Gamma\) is twisted by some higher derivation \(S=(\delta_k)_{k\geq 0}\)
giving rise to commutation laws more general than \eqref{PDOx1} or \eqref{QPDOy1}
(see \cite{FD} for precise definitions). 
The main difficulty lies in this case in the complexity of these commutation laws.

Straightforward calculations show that $z^{-2}=u^{-1}-\frac5{64}\frac{[\xi,\xi]_2}{\xi^6}u-\frac5{64}\frac{[\xi,\xi^2]_3}{\xi^9}u^2+\cdots$ 
In particular \(z^{-2}\not=fu^{-1}+g\) for any \(f,g\in R^\Gamma\), \(f\not=0\).
By uniqueness of the subring \(A^\Gamma\) in \(C^\Gamma\) (see section 5.3 of \cite{FD}),
it follows that the higher derivation \(S\) is actually different of the sequence \(\!(\frac{(2k)\,!}{2^k(k\,!)^2}d^k)_{k\geq 0}\) of \eqref{QPDOy1}
for any derivation \(d\) of \(R^{\Gamma}\).
\end{rem}

Under the assumption of theorem \ref{invCmodu} we can extend part of theorem \ref{*positif} to the space $\mathcal{M}_*$ introduced in remark \ref{MkMl}.

\begin{prop} We suppose that there exists a weight 1 modular form
$\xi$ which is invertible in $R$. Let us define for $k >0$ the vector space morphism $\psi_{-k}: M_{-k} \to C_{-k}^{\Gamma}$ by 
\begin{equation}\label{psimoinsk}
 \psi_{-k}(f)=\psi_{2k}(\xi^{2k})^{-1}\psi_{k}(f\xi^{2k}).
\end{equation}
Then the morphisms $\psi_m$ defined by \eqref{Tpositif} if $m \geq 0$ and by \eqref{psimoinsk} if $m<0$ induce canonically a vector space isomorphism 
$\overline{\Psi}: \mathcal M_* \to C^{\Gamma}$ which defines by transfer
a structure of associative algebra on $\mathcal M_*$.
\end{prop}
\pf Let $k>0$ and let $f \in M_{-k}$, $f \neq 0$. Then $f\xi^{2k} \in M_k$ and $\psi_k(f\xi^{2k}) \in C_k^{\Gamma}$. In the same way $\xi^{2k} \in M_{2k}$, so $\psi_{2k}(\xi^{2k}) \in C_{2k}^{\Gamma}$. Since its $2k$-valuation term is $\xi^{2k}y^{2k}$ and $\xi \in U(R)$ by assumption, we deduce that $\psi_{2k}(\xi^{2k}) \in U(C^{\Gamma})$ and so $\psi_{2k}(\xi^{2k})^{-1} \in C_{-2k}^{\Gamma}$. We have consequently $\psi_{-k}(f) \in C_{-k}^{\Gamma}$, and it is easily checked that its term of valuation $-k$ is $fy^{-k}$. The map $\psi_{-k}$ is clearly linear, and we prove the surjectivity of $\overline {\Psi}$ recursively as in demonstration of theorem \ref{*positif}.
\epf

\begin{rems} 
 \begin{itemize}
  \item[{\rm (i)}] Notice that we don't have in the previous proposition an equivalent of point (ii) of theorem \ref{*positif}, because the morphism $\psi_k$ doesn't satisfy condition (C2) for $k<0$.
  \item[{\rm (ii)}] For $k<0$ an even integer, the map $\psi_k$ defined by \eqref{psimoinsk} can be replaced by the more canonical one introduced previously in formula \eqref{Tnegatifpair}, which satisfies condition (C2). In this case, the isomorphism $\overline{\Psi}_2$ in the remark \ref{defPsi2Barre} is the restriction of $\overline{\Psi}$ to $\mathcal{M}_{*}^{\ev}$.
 \end{itemize}
\end{rems}

\subsection{Additional comment}\label{26}

The action of $\Gamma$ on $B$
and $C$ described in proposition \ref{extensionBCkappa0} and studied throughout the rest of the article 
is based on the choice $r=0$ in formula \eqref{actionsurB}. 
For the same 1-cocycle $s$ defined by \eqref{defdes}, 
we know by example 2 of \ref{examplesr} that another choice
for $r$ could be the map $r':\Gamma\to R$ defined by $r'_{\gamma}=-s_{\gamma}^{-2}d(s_{\gamma}^2)$. Using \eqref{defdes}, we have
$r'_{\gamma} = (cz+d)^{-2}\partial_z((cz+d)^2)
= \frac {2c}{cz+d}$ for any $\gamma = \left(\begin{smallmatrix} a & b \\ c & d\end{smallmatrix}\right) \in \Gamma$. 
In other words, the homographic action of $\Gamma$
on $R$ extends in an action by automorphisms on $B$ defined by:
\begin{equation}\label{action1}
 \trsf{x^{-1}}{}{\gamma}_{1} = (cz+d)^2 (x^{-1} + \frac {2c}{cz+d})
 \qquad\text{ for any }
 \gamma = \left(\begin{smallmatrix} a & b \\ c & d
\end{smallmatrix}\right) \in \Gamma.
\end{equation}

As explained in example 3 of \ref{examplesr}, we can extend for any $\kappa \in \C^*$ the homographic action of $\Gamma$ on $R$ by the action by automorphisms on $B$ defined by:
\begin{equation}\label{actionkappa}
 \trsf{x^{-1}}{}{\gamma}_{\kappa} = (cz+d)^2 (x^{-1} + \kappa \frac {2c}{cz+d})
 \qquad\text{ for any }
 \gamma = \left(\begin{smallmatrix} a & b \\ c & d
\end{smallmatrix}\right) \in \Gamma,
\end{equation}
but the algebraic study of these actions and associated invariant algebras reduces to the case $\kappa=1$. Arithmetical interpretations of these actions have been studied in \cite{CMZ}.
\bibliographystyle{alpha}\nocite{*}
\bibliography{IFPOA}

\end{document}